\documentclass[11pt,reqno]{amsart} 
\usepackage[dvipsnames]{xcolor}
\usepackage{amssymb,
		amsmath,
		amsthm,
		amsfonts,
		enumerate}
\usepackage[
	     breaklinks=true,
		colorlinks=true,
		linkcolor=blue,
		citecolor=ForestGreen,
		urlcolor=blue
		]
{hyperref}
\usepackage{breakurl}\usepackage{soul}
\usepackage{mathtools}
\usepackage[inline]{enumitem}  
\usepackage[margin=1in]{geometry}
\let\oldsquare\square
\usepackage{mathabx}
\usepackage{cite}

\usepackage{subfigure}\usepackage{cases}\usepackage{enumerate}

\makeindex
\newcommand{\n}[1]{\left\Vert#1\right\Vert}
\newcommand{\N}[1]{\left\Vert#1\right\Vert}
\newcommand{\abs}[1]{\left\vert#1\right\vert}
\newcommand{\set}[1]{\left\{#1\right\}}

\newcommand{\R}{\mathbb R}
\newcommand{\rd}{\R^d}
\newcommand{\td}{\mathbb{T}^d}

\newcommand{\Z}{\mathbb Z}

\newcommand{\F}{ {\mathcal F} }

\newcommand{\w}{\widehat{w}_s^{p,q}}

\newcommand{\m}{\mathcal{M}}
\newcommand{\mh}{\widehat{\m}}
\newcommand{\zd}{\mathbb{Z}^d}

\def\TagOnRight

\theoremstyle{plain}
\newtheorem{theorem}{Theorem} [section]
\newtheorem{lemma}[theorem]{Lemma}
\newtheorem{corollary}[theorem]{Corollary}
\newtheorem{proposition}[theorem]{Proposition}

\newtheorem{remark}[theorem]{Remark}

\newtheorem*{A1}{{\bf Theorem A.}}
\newtheorem*{A2}{{\bf Theorem B.}}
\newtheorem*{A3}{{\bf Theorem C}}

\newtheorem*{L1}{{\bf Lemma A}}
\newtheorem*{L2}{{\bf Lemma B}}

\theoremstyle{definition}
\newtheorem{definition}[theorem]{Definition}

\allowdisplaybreaks

\def\({\left(}
\def\){\right)}
\def\<{\left\langle}
\def\>{\right\rangle}

\numberwithin{equation}{section}

\usepackage{tikz}
\usetikzlibrary{patterns}
\begin{document}
\title[Norm inflation  for  fractional Hartree  equations]
{Strong ill-posedness for  fractional Hartree\\  and  cubic NLS Equations} 
\author[D. Bhimani]{Divyang G. Bhimani}
\address{Divyang G. Bhimani
\newline\indent
Department of Mathematics, Indian Institute of Science Education and Research-Pune,\newline\indent Dr. Homi Bhabha Road, Pune 411008, India}
\email{divyang.bhimani@iiserpune.ac.in}

\author[S. Haque] {Saikatul Haque}
\address{Saikatul Haque\newline\indent
Harish-Chandra Research Institute, 
India}
\email{saikatulhaque@hri.res.in}



\maketitle
\begin{abstract} We consider fractional Hartree  and   cubic nonlinear Schr\"odinger  equations on Euclidean space $\mathbb R^d$ and on torus $\td$.  We establish norm inflation (a stronger phenomena than standard  ill-posedness) at every initial data  in   Fourier amalgam spaces with  negative regularity.  In particular, these spaces include Fourier-Lebesgue,  modulation and Sobolev spaces.  We further show that this can be  even worse by exhibiting norm inflation with an infinite loss of regularity. To establish these phenomena, we employ a Fourier analytic approach and introduce new resonant sets corresponding to the fractional dispersion $(-\Delta)^{\alpha/2}$.  In particular,  when  dispersion index $\alpha$ is large enough,  we obtain norm inflation  {above} scaling critical  regularity in some of  these spaces.  
 It turns out that our approach could treat both equations (Hartree and power-type NLS) in a unified manner.  The method   should also  work for a broader range of  nonlinear equations with Hartree-type nonlinearity.
\end{abstract}
\section{Introduction}\label{id}
\subsection{Background}
In this article, we study norm inflation (strong ill-posedness) for the  fractional nonlinear Schr\"odinger equation  of the form
\begin{equation}\label{fhe}
\begin{cases}i\partial_t \psi - (- \Delta)^{\alpha/2} \psi = \mu \(K\ast |\psi|^{2} \)\psi\\
\psi(0, \cdot)=\psi_0
\end{cases} (t,x)\in \mathbb R \times \m
\end{equation}
where $\psi=\psi(t,x)\in \mathbb C, $ and $\mu\in
\{ 1, -1\}.$  Here $\m= \mathbb R^d$ or torus $\mathbb T^d$ and $\widehat{\m}$ denote the Pontryagin dual of $\m,$ i.e. 
\begin{equation*} \widehat{\m}=
\begin{cases}  \rd \quad if \  \m=\rd \\
\mathbb Z^d \quad if \  \m=\td.
\end{cases}
\end{equation*}
The fractional Laplacian is defined as 
\begin{eqnarray*}
\mathcal{F}[(-\Delta)^{\alpha/2}f] (\xi) = (2\pi)^\alpha |\xi|^{\alpha} \mathcal{F}f(\xi) \quad (f\in\mathcal{S}'(\m),  \ \xi \in \widehat{\m}) 
\end{eqnarray*}
where $\mathcal{F}$ denotes the Fourier transform and $\alpha>0.$ Here $K$ is a tempered distribution with 
\begin{equation}\label{hk}
\F K=c_{d,\gamma}|\cdot|^{\gamma-d},  \quad 0<\gamma\leq d.
\end{equation}
We have taken $c_{d,\gamma}$ so that \[K=\begin{cases}
|\cdot|^{-\gamma}&\text{for }\gamma\in(0,d)\\
\delta_0&\text{for }\gamma=d
\end{cases}\]
where $\delta_0$ is  the Dirac distribution with mass at the origin in $\rd.$ Thus, our nonlinearity covers the power-type nonlinearity as well as the Hartree-type nonlinearity. The equation
 \eqref{fhe}  is known as \textit{fractional Hartree equation}  when $0<\gamma<d.$ In particular, for $\alpha=2$ we simply say \textit{Hartree equation}.   And  when  $\gamma=d$  equation  \eqref{fhe} is known as  \textit{cubic   fractional nonlinear Schr\"odinger equation} (fractional cubic NLS).  In particular, for $\alpha=2$ we simply say cubic NLS.

The Hartree equation arises naturally in the dynamics of large quantum systems. It occurs in the context of the mean-field limit of $N$-body dynamics  and developing the theory of stellar collapse,  e.g.,  in the boson particles setting  \cite{frohlich2003mean, lieb1987chandrasekhar}.  \eqref{fhe}  is a fundamental equation of fractional quantum mechanics, a generalization of the standard quantum mechanics extending the Feynman path integral to L\'evy processes \cite{LaskinFra}.  And also,  appears in nonlinear optics, fluids, plasmas; fractional Laplacian  proved to be a powerful tool to describe various physical phenomena \cite{sulem}.  

The equation  \eqref{fhe} enjoys  scaling invariance (dilation symmetry): if $\psi$ solves \eqref{fhe} on $\mathbb R^d$ then $\psi_\lambda$ defined by
$$\psi_{\lambda}(t,x)= \lambda^{\frac{d- \gamma + \alpha}{2}} \psi(\lambda^{\alpha}t, \lambda x)$$ 
  also solves equation \eqref{fhe} with scaled initial data  $\psi_{\lambda}(0,\cdot)=\lambda^{\frac{d-\gamma+ \alpha}{2}} \psi_0(\lambda \cdot).$   The homogeneous Fourier-Lebesgue space 
$\dot{\mathcal{F}L_s^q}(\rd)= \big\{f\in \mathcal{S}'(\mathbb
  R^d): \|f\|_{  \dot{  \mathcal{F}L_s^{q}}}= \| | \cdot |^{s}\F f  \|_{L^{q}}< \infty \big\}$ is invariant under the above scaling when the regularity index $s=s_c$ is given by
\[ s_c= d \left( 1-\frac{1}{q} \right) - \frac{d- \gamma + \alpha}{2}. \]
This  scaling heuristic  suggests that  \eqref{fhe} becomes ill-posed for $s<s_c$.   We note that  \eqref{fhe} may not enjoy  scaling invariance  on the torus,  though  the above $s_c$ still plays a vital role  in our analysis on the torus.     We first recall the notion of well-posedness  in the sense of Hadamard.
\begin{definition}[well-posedness]\label{wpd} Let $X, Y \hookrightarrow \mathcal{S}'(\m)$ be  Banach spaces.  The Cauchy problem \eqref{fhe} is said to be 
 well-posed from  $X$ to $Y$ if, for each bounded subset  $B\subset X$, there exist $T>0$ and a Banach space  $X_{T} \hookrightarrow C([0,T], Y)$ such that
\begin{enumerate}
\item[(a)] for all $\psi_0 \in B$,  \eqref{fhe}  has a unique integral solution \footnote{In this paper,  we say $\psi$ is a solution (or integral solution) to \eqref{fhe} we mean it satisfies \eqref{du}.} $\psi\in X_T$ with  $\psi_{|t=0}=\psi_0$
\item[(b)]  the mapping $\psi_0 \mapsto \psi$ is continuous from $(B, \|\cdot\|_{X})$ to $C([0,T],Y).$
\end{enumerate}
In particular, when $X=Y,$ we say  \eqref{fhe} is locally well-posed in  $X.$  And further,  if $T$ can be chosen arbitrarily large, then we say \eqref{fhe} is globally well-posed in $X$.
\end{definition}
The negation of Definition \ref{wpd} is called a {\textit{lack of
 well-posedness}} or {\textit{instability}}. For instance,  if the solution map $\psi_0 \mapsto \psi$ fails to become continuous, then the problem is not well-posed. This particular instability is referred as  {{\textit{ill-posedness}}}. 

There is an enormous literature on well-posedness  for \eqref{fhe}.  Here we cannot expect to acknowledge  them all.  And thus  we only recall the followings:
\begin{A1}[well-posedness in $\mathcal{F}L^q_s$] \
\begin{enumerate} 
\item  \cite{TohruFra, changxing2008cauchy, HyakunaMulti} The fractional  Hartree equation \eqref{fhe} with $\alpha \in (1,2] $ is locally well-posed in $H^{s}(\rd)$ for $s\geq \frac{\gamma}{2}$ and in $\mathcal{F}L^q (\rd) \ ( q\in ( 2d/(d+\gamma), 2] )$ for $\alpha=2, 0<\gamma<\min (2,d)$ or $\gamma=2, d\geq 3.$ 
\item \cite{SireFrac, Krieger1, GrunrockBi, harrop2020sharp} The fractional cubic NLS \eqref{fhe} with $\alpha \in (0,2)$ except for $\alpha=1$  is locally well-posed in $H^s(\rd)$ for $s>s_c$ for $d\geq 2$ and in $H^s(\mathbb R) \ (s\geq 1/2)$ for $\alpha=1$  and in $\mathcal{F}L^q(\R) \ (1\leq q < \infty)$ for $\alpha=2$ and $\mu=-1.$  The 1D cubic NLS is globally well-posed in $H^{s}(\mathbb R)$ for $s>-\frac{1}{2}$.
\end{enumerate}   
\end{A1} 

\begin{A2}[well-posedness in modulation spaces $M^{p,q}_s$]\label{wmt} Let $1/2<\alpha \leq 2.$
\begin{enumerate}
\item \cite{BhimaniFrac, bhimani2020hartree}  The fractional Hartree equation \eqref{fhe}  is locally well-posed in $M^{p,q}(\rd) \ (1\leq p \leq 2, 1\leq q \leq 2d/(d+\gamma))$ and, in  $M^{2,q}(\rd)$  with extended  $ \alpha \in (2, \infty)$. 
\item \cite{bhimani2016functions, baoxiang2006isometric, BenyiLocal}. The fractional  cubic NLS \eqref{fhe} is locally well-posed in  $M^{p,1}_s(\rd) \ (1\leq p \leq \infty, s\geq 0)$ and in $M^{p,q}_s(\rd) \ (1\leq p, q \leq \infty, s> d(1-1/q)).$  Furthermore, the range of $\alpha$ can be extended to $(2, \infty)$ if $p=2$.
\item  \cite{OhWang2020, wang2007global} The 1D cubic NLS \eqref{fhe} is globally well-posed in $M^{2,q}(\mathbb R)$ for $2\leq q < \infty.$ The cubic NLS \eqref{fhe} is globally well-posed (for sufficiently small data) in $M_s^{2,1}(\rd)$ for $d\geq 2$ and $s\geq 0.$ 
\end{enumerate}
\end{A2}  

\begin{A3}[well and ill-posedness on torus]\label{torus}\
\begin{enumerate}
\item \cite{Bourgain} The cubic NLS  \eqref{fhe} is globally well-posed in $L^2(\mathbb T).$
 \item \cite{Molinet, Christ, Burq, kishimoto2019remark, oh2017remark}  The cubic NLS  \eqref{fhe}  is ill-posed in $H^{s}(\mathbb T)$ for $s<0.$
 \item  \cite{hirocircle} The cubic fractional  NLS \eqref{fhe} is strongly  ill-posed (norm inflation) in $H^{s}(\mathbb T)$ if $s<s_c$, $s\neq0$, $0<\alpha <2$ or  if $s\leq s_c$, $\alpha\geq2$. 
\item \cite{carles2017norm} The cubic NLS \eqref{fhe} is strongly ill-posed (with infinite loss of regularity) in $\mathcal{F}L^p_s(\mathbb T)$ for $s<-2/3$ and in $\mathcal{F}L^p_s(\mathbb T^d)$ for $s<0$ and $d\geq 2;$ for $p\in [1, \infty].$
\end{enumerate}
\end{A3}

The purpose of the present paper is to complement these positive results by establishing  ill-posedness for \eqref{fhe} and also  to  improve and extend  several  negative (i.e. ill-posedness)  results  from \cite{oana,carles2012geometric,hirocircle, kishimoto2019remark, bhimani2020norm, changxing2008cauchy, han2019global}.  
Now,  in order to state and discuss our main  results (in Subsection \ref{mr}),  we  shall first  briefly  introduce  the following function spaces. 
\subsection{Function Spaces}
Recently in \cite{JH},     Oh and Forlano have introduced Fourier amalgam spaces $\widehat{w}^{p,q}_s(\m)$  (with $1\leq p, q \leq \infty,  s\in \mathbb R)$ 
$$\widehat{w}^{p,q}_s(\m)=\left \{ f\in \mathcal{S}'(\m): \|f\|_{\widehat{w}^{p,q}_s}=   \left\| \left\lVert \chi_{n+Q}(\xi)\F f(\xi)\right\rVert_{L_\xi^p(\mh)}  \langle n \rangle^s \right\|_{\ell^q_n(\zd)}< \infty \right\},
$$
where $Q=(-\frac{1}{2},\frac{1}{2}]^d$   and  $\langle \cdot \rangle^s= (1+|\cdot|^2)^{s/2}.$  Their  motivation was to study  well-posedness of 1D cubic NLS in these spaces, \cite{forlano2020deterministic}. We note that these spaces are  Fourier image of Wiener amalgam spaces,  see Remark \ref{chr}.  And as far as we are aware,  there are no well-posedness results for Hartree equation in  these spaces  so far except the case $p=q=2$ when it coincide with $H^s(\m)$.
On the  other hand,  the modulation $M^{p,q}_s(\m)$ 
 spaces were introduced by Feichtinger in early  1980's  in \cite{feichtinger1983modulation} and proved to be very fruitful to study  nonlinear dispersive equations in the last decades \cite{wang2011harmonic, wang2007global, baoxiang2006isometric, OhWang2020, bhimani2016functions}.
To recall their definitions, let   $\rho \in \mathcal{S}(\R^d),$  $\rho: \R^d \to [0,1]$  be  a smooth function satisfying   $\rho(\xi)= 1$ if $|\xi|_{\infty}=\max(|\xi_1|,...,|\xi_d|)\leq \frac{1}{2} $ and $\rho(\xi)=0$ if $|\xi|_{\infty}\geq  1.$ Let  $\rho_n$ be a translation of $\rho,$ that is,
$ \rho_n(\xi)= \rho(\xi -n), n \in \Z^d$
and denote 
$\sigma_{n}(\xi)=
  \frac{\rho_{n}(\xi)}{\sum_{\ell\in\Z^{d}}\rho_{\ell}(\xi)},  n
  \in \Z^d.$
Then the frequency-uniform decomposition operators can be defined by 
\[\oldsquare_n = \mathcal{F}^{-1} \sigma_n \mathcal{F}, \quad n\in\Z^d. \]
Now 
the modulation  space $M^{p,q}_s(\m)$ are defined by the following norm
\begin{equation*}
\|f\|_{M^{p,q}_s(\m)}=   \left\| \left\lVert
  \oldsquare_nf\right\rVert_{L_x^p(\m)} \langle n \rangle ^{s} \right\|_{\ell^q_n(\zd)}. 
\end{equation*}
 See Remark \ref{stft} for equivalent characterization via short-time Fourier transform (STFT) of these spaces.  The Fourier-Lebesgue spaces $\mathcal{F}L^q_s(\m)$ is defined by 
$$\mathcal{F}L^q_s(\m)=\left\{f\in \mathcal{S}'(\m): \mathcal{F}f \  \langle  \cdot  \rangle^s  \in L^q (\widehat{\m})\right \}.$$ It is interesting to note that $\widehat{w}^{p,q}_s(\m)$ recapture several known spaces:
\begin{equation*} 
   \widehat{w}_s^{p,q}(\m)= \begin{cases}
   \mathcal{F}L_s^q(\m) \ (\text{Fourier-Lebesgue space}) &\text{if } p=q \ \\
   M^{2,q}_s(\m) \ (\text{modulation space}) &\text{if }  p=2 \ \\
   H^{s}(\m) \ (\text{Sobolev space})   &\text{if }  p=q=2 \  \\
   \mathcal{F}L_s^q(\m)=M_s^{p,q}(\m) &\text{if }    \m =\td.
    \end{cases}
    \end{equation*}
\subsection{Main Results}{\label{mr}}
\begin{theorem}\label{mt1} Let $0< \gamma \leq d, p,q \in \mathbb [1, \infty]$ and $\alpha>0.$   Assume that  $s$ satisfies  one of the following conditions (see Remark \ref{ur}):
 \begin{enumerate}[label=\fbox{\arabic*}]
 \item $ s<\min\(-2(d-\gamma),s_c+
\max\big(\frac{(d-\gamma)}{2}(1+\frac{2d}{sq'}),(d-\gamma)(1+\frac{2d}{sq'})\big)\)$\label{mt11}.
 \item $s< \min\(-2(d-\gamma),-\frac{d-\gamma+\alpha}{3}\).$\label{mt12}
 \item  $s<-\frac{d}{q}-(d-\gamma)$ and $s<-\frac{\alpha+d-\gamma}{3}+\begin{cases}\frac{d}{2}(\frac{2}{3}-\frac{1}{q}) \text{ if } q\leq\frac{3}{2}\\
 d(\frac{2}{3}-\frac{1}{q}) \text{ if } q\geq\frac{3}{2}
 \end{cases}$.\label{mt13}
 \item  $s=-\frac{d}{q}\ \text{and}\ d+2(d-\gamma)<\alpha<2d+2(q-1)(\gamma-d).$
\label{mt14}
 \end{enumerate}
Then norm inflation   occurs to \eqref{fhe} everywhere in $\widehat{w}_s^{p,q}(\m)$,  that is,
for any  base $\psi_0\in \widehat{w}_s^{p,q}(\m)$ and $\varepsilon>0$ there exist a smooth $\psi_{0,\varepsilon} \in \widehat{w}_s^{p,q}(\m)$ and  $T>0$ satisfying  
 \[ \|\psi_0-\psi_{0,\varepsilon}\|_{\widehat{w}_{s}^{p,q}}< \varepsilon, \quad 0<T< \varepsilon\] 
 such that the corresponding smooth solution $\psi_\varepsilon$ to $\eqref{fhe}$ with data $\psi_{0,\varepsilon}$ exists on $[0,T]$ and 
 \[ \|\psi_\varepsilon(T)\|_{\widehat{w}^{p,q}_s}> \frac{1}{\varepsilon}
 .\] 
Moreover, the  hypothesis in  fourth scenario \ref{mt14} can be relaxed to 
\begin{enumerate}[resume]
\item[\fbox{4'}] $s=-\frac{d}{q}\quad\text{and}\quad
   0<\alpha<\begin{cases}
   \frac{3d}{q}&\text{for }q\leq\frac{3}{2}\\ 
2d&\text{for }q\geq\frac{3}{2}
   \end{cases}$
(with $q<\infty$,  the case $\alpha=2d$ can also  be included)\end{enumerate}  to get norm inflation for fractional  cubic NLS.
\end{theorem} 
\begin{figure}
\subfigure[The case $d=1$, $\gamma=\frac{1}{2}$, $q=1$] 
{ 
\begin{tikzpicture}[scale=.75]
\fill [gray!50!white] (0,-1) -- plot[domain=2.5:7.5]  (\x, {-1/6-\x/3})--(7.5,-3.25)--(0,-3.25);
\fill [gray!50!white] (0,-1) -- plot[domain=3:7.5]  (\x, {1/2-\x/2})--(7.5,-3.25)--(0,-3.25);
\draw[][->] (-.25,0)--(9,0) node[anchor=south] {\tiny{$\alpha$}};
\draw[][<-] (0,.25)--(0,-3.5) node[anchor=east] {\tiny{$s$}};
\draw[dotted][-] (0,-1)--(4,-1) node[anchor=west] {\tiny{$s=-2(d-\gamma)=-1$}};
\draw[dotted][-] (0,-1.5)--(4.8,-1.5) node[anchor=west] {\tiny{$s=-\frac{d}{q}-(d-\gamma)=-\frac{3}{2}$}};
\draw[dotted][->] [domain=.0:7] plot ({\x},{-.5*\x});
\draw[dashed][->] [domain=1:7.5] plot ({\x},{-\x/2+1/2});
\draw[dashed][->] [domain=0:7.5] plot ({\x}, {-1/6-\x/3})node[anchor=south] {\tiny{\fbox{2}}};
\draw[dotted][->] [domain=0:7.5] plot ({\x}, {-1/6-\x/3-1/6})node[anchor=west] {\tiny{\textcolor{gray}{\fbox{3a}}}};
\draw[<-](2,-.5)--(2.5,-.5) node[anchor=west]{\tiny{\fbox{1$b$}}};
\draw(2.5,-2.5) node[anchor=west]{\tiny\textcolor{gray}{{\fbox{1$a$}}}};\draw[->](3.5,-2.5)--(5.0,-2.5);
\end{tikzpicture}\label{f1a}
} 
\subfigure[The case $d=1$, $\gamma=\frac{9}{10}$, $q=2$]  
{ 
\begin{tikzpicture}[scale=1.5]
\fill [gray!50!white] (0,-3/5)-- plot[domain=2.2:7.5/2] (\x,{-1/30-\x/3+1/6})--(7.5/2,-3.25/2)--(0,-3.25/2);
\fill [gray!50!white] (0,-.2) -- plot[domain=.9:7.5/2]  (\x, {.5*(-\x/2+1/2)-.5*((-\x/2+1/2)^2+4/20)^(.5)})--(0,-3.25/2);
\draw[][->] (-.5/4,0)--(4.5,0) node[anchor=south] {\tiny{$\alpha$}};
\draw[][<-] (0,.125)--(0,-3.5/2) node[anchor=east] {\tiny{$s$}};
\draw[dotted][-] (0,-1/5)--(2.5,-1/5) node[anchor=west] {\tiny{$s=-2(d-\gamma)=-\frac{1}{5}$}};
\draw[dotted][-] (0,-1/2-1/10)--(2.4,-1/2-1/10) node[anchor=west] {\tiny{$s=-\frac{d}{q}-(d-\gamma)=-\frac{3}{5}$}};
\draw[dashed][->] [domain=.0:7.5/2] plot ({\x},{.5*(-\x/2+1/2)-.5*((-\x/2+1/2)^2+4/20)^(.5)});
\draw[dotted][->] [domain=.0:4] plot ({\x},{.5*(-\x/2+11/20)-.5*((\x/2-11/20)^(2)+8/20)^(.5)})node[anchor=west] {\tiny{\textcolor{gray}{\fbox{1$b$} }}};
\draw[dotted][->] [domain=0:7.5/2] plot ({\x}, {-1/30-\x/3})node[anchor=west] {\tiny{\textcolor{gray}{\fbox{2}}}};
\draw[dashed][->] [domain=0.4:7.5/2] plot ({\x}, {-1/30-\x/3+1/6});
\draw[<-](1.8,-.5)--(2.3,-.4) node[anchor=west]{\tiny{\fbox{1$a$}}};
\draw[<-](3.4,-1)--(3.8,-1) node[anchor=west]{\tiny{\fbox{3$b$}}};
\end{tikzpicture}\label{f1b}
}

\begin{align*}\tiny{
\fbox{1a} }& \ \tiny{\text{$s=s_c+\frac{1}{2}(d-\gamma)(1+\frac{2d}{sq'})$}}&\tiny{\fbox{1b}}&\ \tiny{\text{$s=s_c+(d-\gamma)(1+\frac{2d}{sq'})$}}&\tiny{\fbox{2} }&\ \tiny{\text{$s=-\frac{d-\gamma+\alpha}{3}$}}\\
\tiny{\fbox{3a} }&\ \tiny{\text{$s=-\frac{d-\gamma+\alpha}{3}+\frac{d}{2}(\frac{2}{3}-\frac{1}{q})$}}&\tiny{\fbox{3b} }&\ \tiny{\text{$s=-\frac{d-\gamma+\alpha}{3}+d(\frac{2}{3}-\frac{1}{q})$}}
\end{align*}
\caption{ Norm inflation (NI for short)  in Hartree case occurs in the shaded regin (the unbounded region between $s$-axis and the curves, lines). }\label{f1}
\end{figure}

\begin{figure}
\subfigure[The case $1\leq q<\frac{3}{2}$ (figure is exact for $q=\frac{5}{4}$)]  
{ 
\begin{tikzpicture}
\fill [gray!50!white] (0,0) -- plot[domain=2/5:7.5/2]  (\x, {1/5-\x/2})--(7.5/2,-3.25/2)--(0,-3.25/2);
\fill[] [gray!90!white] (0,0) -- plot[domain=2/5:7.5/2]  (\x, {-\x/3})--(7.5/2,-3.25/2)--(0,-3.25/2);
\draw[][->] (-.25,0)--(6,0) node[anchor=south] {\tiny{$\alpha$}};
\draw[][<-] (0,.25)--(0,-2) node[anchor=east] {\tiny{$s$}};
\draw[dashed](0,-4/5)--(4,-4/5)node[anchor=west] {\tiny{{$s=-\frac{d}{q}$}}};
\draw[dashed](2/5,-1.7)--(2/5,0) node[anchor=south] {\tiny{$\frac{2d}{q'}$}};
\draw[dashed](6/5,-1.7)--(6/5,0) node[anchor=south] {\tiny{$\frac{6d}{q'}$}};
\draw[dashed][->] [domain=2/5:4] plot ({\x}, {1/5-\x/2})node[anchor=west] {\tiny{{$s=s_c$}}};
\draw[dashed][->] [domain=0:4] plot ({\x}, {-\x/3})node[anchor=west] {\tiny{{$s=-\frac{\alpha}{3}$}}};
\end{tikzpicture}
}
\subfigure[The case $ q>\frac{3}{2}$ (figure is exact for $q=2$)]  
{ 
\begin{tikzpicture}
\fill [gray!50!white] (0,0) -- plot[domain=1:7.5/2]  (\x, {1/2-\x/2})--(7.5/2,-3.25/2)--(0,-3.25/2);
\fill [gray!50!white] (0,0) -- plot[domain=2/5:7.5/2]  (\x, {1/6-\x/3})--(7.5/2,-3.25/2)--(0,-3.25/2);
\fill[] [gray!90!white] (0,0) -- plot[domain=2/5:7.5/2]  (\x, {-\x/3})--(7.5/2,-3.25/2)--(0,-3.25/2);
\draw[][->] (-.25,0)--(6,0) node[anchor=south] {\tiny{$\alpha$}};
\draw[][<-] (0,.25)--(0,-2) node[anchor=east] {\tiny{$s$}};
\draw[dashed](0,-1/2)--(4,-1/2)node[anchor=west] {\tiny{{$s=-\frac{d}{q}$}}};
\draw[dashed](1,-1.7)--(1,0) node[anchor=south] {\tiny{$\frac{2d}{q'}$}};
\draw[dashed](2,-1.7)--(2,0) node[anchor=south] {\tiny{$2d$}};
\draw[dashed][->] [domain=1:4] plot ({\x}, {1/2-\x/2})node[anchor=west] {\tiny{{$s=s_c$}}};
\draw[dotted][->] [domain=0:4] plot ({\x}, {-\x/3})node[anchor=west] {\tiny{\textcolor{gray}{$s=\ \ \ -\frac{\alpha}{3}$}}};
\draw[dashed][->] [domain=0.5:4] plot ({\x}, {1/6-\x/3})node[anchor=west] {\tiny{{$s=d(\frac{2}{3}-\frac{1}{q})-\frac{\alpha}{3}$}}};
\end{tikzpicture}
}
\caption{\small{NI for power-type ($d=\gamma=1$) NLS occurs in the shaded (black \& gray) region, 
 the point $(2d,-\frac{d}{q})$ is included for $1\leq q<\infty$ (Theorem \ref{mt1}, Corollary \ref{c1}).   If $\alpha>\frac{6d}{q'}$ with $1\leq q<\frac{3}{2}$ or $\alpha>2d$ with $q\geq\frac{3}{2}$ we have achieved NI in $\w$ with some $s>s_c.$ NI with infinite loss of regularity occurs in the darked region.
 }
 }\label{f2}
\end{figure}
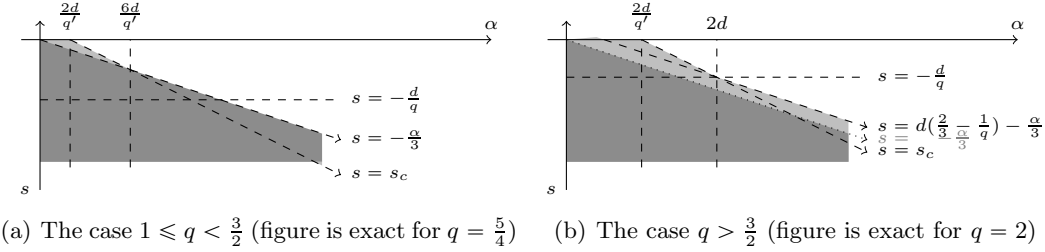

\begin{remark}\label{ur}
Each condition \ref{mt11}-\ref{mt14} in Theorem \ref{mt1} has its own advantage with different choices of parameters $d,\gamma,q, \alpha,s$. This is demonstrated via Figures \ref{f1} and \ref{f2}. In Figure $\ref{f1a}$,  Cases \ref{mt11} (with $\fbox{1b}$) and \ref{mt12}  cover more range of $s$ with certain choices of $\alpha$, whereas  in Figure $\ref{f1b}$  Cases  \ref{mt11} (with $\fbox{1a}$), \ref{mt13}  has similar advantage.  Figure \ref{f2} describes  case $\gamma=d$, i.e.   for  cubic fractional NLS \eqref{fhe}.
\end{remark}

Now we discuss several interesting features of Theorem \ref{mt1}:
\begin{itemize}
\item In spite of much progress for NLS on a torus (for e.g. Theorem C), there are no ill-posedness results for Hartree equation on torus as far as we are aware. Theorem \ref{mt1} thus  is the first ill-posedness result  for the Hartree equation \eqref{fhe} on the torus $\td$.  See,  for example,   \cite{Nakamura, Vedran}  for work on the Hartree equation on torus.
\item In   \cite[Theorem 1.5]{oana} (the final case of the theorem), Choffrut and  Pocovnicu  proved  NI for  fractional 1D cubic NLS \eqref{fhe} based at 0 in $H^s(\mathbb R)$ for $s< \frac{1-2\alpha}{6}$ and $\alpha>2.$ Note that 
the particular case ($q=2, d=\gamma=1$) of Theorem \ref{mt1} \ref{mt13} improves this result in the sense that we have achieved NI at general initial data $\psi_0$. 

\item For $\alpha>2,$ one has $s_c<-\frac{\alpha}{3}+(\frac{2}{3}-\frac{1}{2})$.
Therefore, Theorem \ref{mt1} \ref{mt13} extends the result of Oh-Wang \cite[Theorem 1.2]{hirocircle}  
 mentioned in Theorem C.

\item Bhimani-Carles  \cite[Theorems 1.2, 1.5]{bhimani2020norm} proved norm inflation  based at zero initial data for NLS  \eqref{fhe}  in $\mathcal{F}L^q_s(\rd)$ and in $M_s^{2,q}(\rd)$  for $s<\min (0, s_c).$  In dimension one,  Theorem \ref{mt1} $\fbox{4'}$ improves this result by reaching up to critical  regularity $s_c=-1/q$.  And also  extends  to  fractional NLS \eqref{fhe}. 

\item By Theorem B (1)  fractional  Hartree equation  \eqref{fhe} is well posed in $M_0^{2,1}(\rd)$. Theorem \ref{mt1} leaves a gap at least  for  
$s\in \big(\max(-2(d-\gamma),-(d-\gamma)-\frac{d}{q}), 0\big)$ ($-\frac{d}{q}$ is excluded in some cases).  

\item By Theorem B (3), 1D cubic NLS is globally well-posed $M^{2,q}(\mathbb R)$ for $q\geq2$. Thus Theorem \ref{mt1}  leaves a gap 
$\big[\max\{s_c,-\frac{\alpha}{3},\min\{-\frac{\alpha}{3}+\frac{2}{3}-\frac{1}{q},-\frac{1}{q}\}\},0\big)$ ($-\frac{1}{q}$ is excluded in some cases). However, in dimension one, a weaker version of ill-posedness is  recently proved for such $s$'s: the solution map is not $C^3$ in $M_s^{p,q}(\R)$ for all $s<0$, see \cite[Theorem 6.1]{klaus2022wellposedness}.

\end{itemize}
\begin{corollary} \label{c1} Let $\gamma= d$ and assume that $\alpha$ and $s$ satisfy the  following  (see the shaded region in Figure \ref{f2}):
\[
s<\begin{cases}
0&\text{if }\ 0<\alpha\leq\frac{2d}{q'}\\
s_c&\text{if }\frac{2d}{q'}\leq\alpha\leq\frac{6d}{q'}\\
-\frac{\alpha}{3}&\text{if }\frac{6d}{q'}\leq\alpha<\infty
\end{cases}\ \text{for }q\in[1,\frac{3}{2}]\text{ and }
s<\begin{cases}
0&\text{if }\ 0<\alpha\leq\frac{2d}{q'}\\
s_c&\text{if }\frac{2d}{q'}\leq\alpha\leq2d\\
-\frac{\alpha}{3}+d(\frac{2}{3}-\frac{1}{q})&\text{if }2d\leq\alpha<\infty
\end{cases}\ \text{for }q\in[\frac{3}{2},\infty].
\]
Then norm inflation for \eqref{fhe} occurs everywhere in $\w(\m).$ Moreover, with $q\in[1,\infty)$ and $\alpha=2d$, norm inflation for \eqref{fhe} occurs everywhere in $\widehat{w}_{s_c}^{p,q}(\m)$.
\end{corollary}
 As a consequence of Theorem \ref{mt1}  we the following instability for \eqref{fhe}. 
 \begin{corollary}\label{cr} Under the assumption of Theorem \ref{mt1},
  for any $T>0,$ the solution map $\widehat{w}_s^{p,q}(\m)  \ni \psi_0\mapsto \psi \in C([0, T],  \widehat{w}_{s}^{p,q}(\m))$ for  \eqref{fhe}  is discontinuous everywhere  in $\widehat{w}^{p,q}_{s}(\m)$.
 \end{corollary}
 
Corollary \ref{cr} states that \eqref{fhe} experiences ill-posedness.  Now,  we shall show that the situation  is even worse in the sense that  the solution map fails to be continuous  from  $\widehat{w}^{p,q}_s(\m)$ to $\widehat{w}^{p,q}_{\theta}(\m)$ \textit{even for (some) $\theta\neq s$}.  And so, a loss of regularity is {exhibited}.  This phenomenon  is  initially exhibited   by  Lebeau  \cite{MR2145023} in the case of the wave equation and  then by Carles in \cite{carles2007geometric}, Thomann in \cite{thomann2008instabilities} and Alazard-Carles in  
 \cite{alazard2009loss}
 for NLS \eqref{fhe}.     In  all these works,  loss of regularity (finite or infinite) was established based at zero initial data with some positive regularity.   In the next theorem,  we establish norm inflation with infinite loss of regularity based at general initial data with negative regularity.
\begin{theorem}\label{mt2} Under the assumption \ref{mt12} in Theorem \ref{mt1}, 
norm inflation with infinite loss of regularity   occurs to \eqref{fhe} everywhere in $\widehat{w}_s^{p,q}(\m)$,  that is,
for any $\psi_0\in \widehat{w}_s^{p,q}(\m)$ and $\varepsilon>0$ there exist a smooth $\psi_{0,\varepsilon} \in \widehat{w}_s^{p,q}(\m)$ and  $T>0$ satisfying  
 \[ \|\psi_0-\psi_{0,\varepsilon}\|_{\widehat{w}_{s}^{p,q}}< \varepsilon, \quad 0<T< \varepsilon\] 
 such that the corresponding smooth solution $\psi_\varepsilon$ to $\eqref{fhe}$ with data $\psi_{0,\varepsilon}$ exists on $[0,T]$ and 
 \[ \|\psi_\varepsilon(T)\|_{\widehat{w}^{p,q}_{\theta}}> \frac{1}{\varepsilon} \quad \text{for all} \ \theta \in \mathbb R
 .\] 
In particular,  for any $T>0,$ the solution map $\widehat{w}_s^{p,q}(\m)  \ni \psi_0\mapsto \psi \in C([0, T],  \widehat{w}_{\theta}^{p,q}(\m))$ for  \eqref{fhe}  is discontinuous everywhere  in $\widehat{w}^{p,q}_{s}(\m)$
\end{theorem} 
Theorem \ref{mt2} deserve several comments:
\begin{itemize}

\item In \cite{carles2017norm},  Carles and  Kappeler  proved infinite loss of regularity  (via geometric approach) for 1D cubic NLS  based at zero initial data  in $\mathcal{F}L^q(\mathbb T)$ for $s<-2/3.$ (See Theorem C(4)). Theorem \ref{mt2} strengthens this   by establishing an infinite loss of regularity at general initial data.  And also extends to fractional NLS.

 \item In \cite[Theorem 1.6]{bhimani2020norm}, the first author with Carles established norm inflation  with infinite loss of regularity for cubic NLS \eqref{fhe} in  $M^{p,q}_s(\R^d) \ (d\geq 2, 1\leq p, q \leq \infty, s<-1/3)$ via geometric optics approach. Theorem \ref{mt2} establish norm inflation with infinite loss of regularity  for fractional cubic NLS \eqref{fhe} in $\w$ for any $d\geq 1$ via Fourier analytic approach. However, we note that  the   restriction on the first index $p=2$ in the case of modulation spaces  in Theorem \ref{mt2} is due to our approach and we  believe that  this restriction should be relaxed by taking any $1\leq p \leq \infty.$

\item To the best of the authors' knowledge,   infinite  loss of regularity  has  not been  previously studied for Hartree equation \eqref{fhe} even on 
$\mathbb R^d.$ Thus, Theorem \ref{mt2} is the first result in this direction. 


\end{itemize}

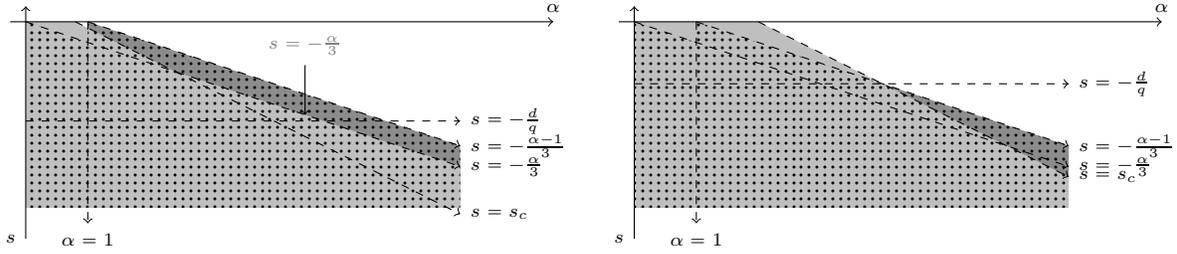
\begin{figure}
\subfigure[The particular case $q=\frac{5}{4}$]  
{ 
\begin{tikzpicture}[scale=.825]
\fill[] [gray!90!white] (1,0) -- plot[domain=1:7]  (\x, {-\x/3+1/3})--(7,-3)--(1,-3);
\fill [gray!50!white] (0,0) -- plot[domain=4/5:7.5/2]  (\x, {2/5-\x/2})--(7.5/2,-3.25/2)--(0,-3.25/2);
\fill[] [gray!50!white] (0,0) -- plot[domain=0:7]  (\x, {-\x/3})--(7,-3)--(0,-3);
\fill[pattern=dots]  (0,0) -- plot[domain=0:1]  (\x, {-\x/3})--(1,-3)--(0,-3);
\fill[pattern=dots]  (1,0) -- plot[domain=1:7]  (\x, {-\x/3+1/3})--(7,-3)--(1,-3);
\draw[][->] (-.25,0)--(8.5,0) node[anchor=south] {\tiny{$\alpha$}};
\draw[][<-] (0,.25)--(0,-3.5) node[anchor=east] {\tiny{$s$}};
\draw[dashed][->] [domain=0:-3.25] plot ({1}, {\x})node[anchor=north] {\tiny{{$\alpha=1$}}};
\draw[dashed][->] [domain=4/5:7] plot ({\x}, {2/5-\x/2})node[anchor=west] {\tiny{{$s=s_c$}}};
\draw[dashed][->] [domain=1:7] plot ({\x}, {-\x/3+1/3})node[anchor=west] {\tiny{{$s=-\frac{\alpha-1}{3}$}}};
\draw[dashed][->] [domain=0:7] plot ({\x}, {-\x/3})node[anchor=west] {\tiny{{{$s=-\frac{\alpha}{3}$}}}};
\draw[dashed][->] [domain=0:7] plot ({\x}, {-8/5})node[anchor=west] {\tiny{{$s=-\frac{d}{q}$}}};
\draw[<-](4.5,-1.5)--(4.5,-.7)node[anchor=south]{\tiny{{\textcolor{gray}{$s=-\frac{\alpha}{3}$}}}};
\end{tikzpicture}
}
\subfigure[The particular case $q=2$]  
{ 
\begin{tikzpicture}[scale=.825]
\fill [gray!90!white] (4,-1) -- plot[domain=5:7]  (\x, {-\x/3+1/3})--(7,-3)--(0,-3);
\fill [gray!50!white] (0,0) -- plot[domain=2:7]  (\x, {1-\x/2})--(7,-3)--(0,-3);
\fill[pattern=dots]  (0,0) -- plot[domain=0:1]  (\x, {-\x/3})--(1,-3)--(0,-3);
\fill[pattern=dots]  (1,0) -- plot[domain=1:7]  (\x, {-\x/3+1/3})--(7,-3)--(1,-3);
\draw[][->] (-.25,0)--(8.5,0) node[anchor=south] {\tiny{$\alpha$}};
\draw[][<-] (0,.25)--(0,-3.5) node[anchor=east] {\tiny{$s$}};
\draw[dashed][->] [domain=0:-3.25] plot ({1}, {\x})node[anchor=north] {\tiny{{$\alpha=1$}}};
\draw[dashed][->] [domain=2:7] plot ({\x}, {2/2-\x/2})node[anchor=west] {\tiny{{$s=s_c$}}};
\draw[dashed][->] [domain=1:7] plot ({\x}, {-\x/3+1/3})node[anchor=west] {\tiny{{$s=-\frac{\alpha-1}{3}$}}};
\draw[dashed][->] [domain=0:7] plot ({\x}, {-\x/3})node[anchor=west] {\tiny{{$s=-\frac{\alpha}{3}$}}};
\draw[dashed][->] [domain=0:7] plot ({\x}, {-2/2})node[anchor=west] {\tiny{{$s=-\frac{d}{q}$}}};
\end{tikzpicture}
}

\caption{\small{$\gamma=d=2$: NI occurs in the shaded region (the point $(2d,-\frac{d}{q})$), improvement by Theorem \ref{mt3} is indicated by darker region. NI with infinite loss of regularity occurs in region indicated by dotted region.
 }
 }\label{f4}
\end{figure}

With the existence of resonance set (see \eqref{rv}) we have the following improvement to Theorems \ref{mt1} and \ref{mt2} for certain choices of $d$ and $\alpha$:
\begin{theorem}[Improvement]\label{mt3}
Assume $d=\alpha=1$ or $d\geq2,\alpha\geq1$.  Let $0< \gamma \leq d$ and $p,q \in \mathbb [1, \infty].$
\begin{enumerate}
\item\label{mt31} Assume that $s$ satisfies  one of the following conditions:
\begin{enumerate}[label=\fbox{\alph*}]
\item $s<\min\(-2(d-\gamma),s_c+\max\big(\frac{(d-\gamma)}{2}(1+\frac{2d}{sq'}-\frac{1}{s}),(d-\gamma)(1+\frac{2d}{sq'}-\frac{1}{s})\big)\).$\label{mt3a}
\item $s<\min\(-2(d-\gamma),-\frac{d-\gamma+\alpha-1}{3}\).$\label{mt3b}
\item $s<-\frac{d}{q}-(d-\gamma)$ and $s<-\frac{d-\gamma+\alpha}{3}+\frac{1}{6}+\frac{d}{2}(\frac{2}{3}-\frac{1}{q}) \text{ if } q\leq\frac{3d}{2d-1}$.\label{mt3c}
\item $s=-\frac{d}{q},\ 1+d+2(d-\gamma)<\alpha\leq 2d+[2-(2-\frac{1}{d})q](d-\gamma).$
\label{mt3d}
\end{enumerate}
Then norm inflation occurs to \eqref{fhe} everywhere in $\widehat{w}_s^{p,q}(\m)$.
Moreover,  for fractional cubic NLS, the hypothesis of  \ref{mt3d} can be relaxed to 
\begin{enumerate}[resume]
\item[\fbox{d'}]$s=-\frac{d}{q},\ 1<\alpha<
1+ \frac{3d}{q}\text{ for }q\leq\frac{3d}{2d-1}
$
(the case $\alpha=2d$ can also be  included by excluding $q=\infty$), to have norm inflation everywhere in $\widehat{w}_s^{p,q}(\m)$.\end{enumerate}
\item\label{mt32} Under the assumption \ref{mt3b} in part \eqref{mt31}
norm inflation with infinite loss of regularity   occurs to \eqref{fhe} everywhere in $\widehat{w}_s^{p,q}(\m)$.
\end{enumerate}
\end{theorem}
Theorem \ref{mt3} deserve several comments:
\begin{itemize}
\item Since $s<0$, it is easy to see that part \ref{mt3a} of Theorem \ref{mt3} improves (enlarges the range of $s$ where NI occurs) part  \ref{mt11} of Theorem \ref{mt1}.  This is also reflected in the Figure \ref{f3}. On the other hand,  part  \ref{mt3b} improves part \ref{mt12} of Theorem \ref{mt1}.  Part  \ref{mt3c} improves  part \ref{mt13} of Theorem \ref{mt1}  if $1\leq q\leq\frac{3}{2}$. Compared to Theorem \ref{mt1} \ref{mt14} \fbox{4'}, the conditions \ref{mt3d} and \fbox{d'} of Theorem \ref{mt3} improve the range of $\alpha$. Finally Theorem \ref{mt3}\eqref{mt32} improves Theorem \ref{mt2}.

\item In Figure \ref{f4}, using Theorems \ref{mt1}, \ref{mt2} and \ref{mt3} for the case $\gamma=d=2$, i.e. two dimensional power-type NLS case, we have indicated where NI occurs, where NI with infinite loss of regularity occurs and where we have used the improvement through Theorem \ref{mt3}.   
\item Miao-Xu-Zhang \cite[Theorem 4.2]{changxing2008cauchy} established 
 NI for  Hartree equation \eqref{fhe} in $H^s(\rd)$ for $s\leq-\frac{d}{2}$.  
 Theorem \ref{mt3} with the assumption \ref{mt3b} extends this result whenever $d+2\gamma\geq2$, $\frac{3}{4}d<\gamma< d$.  

\end{itemize}
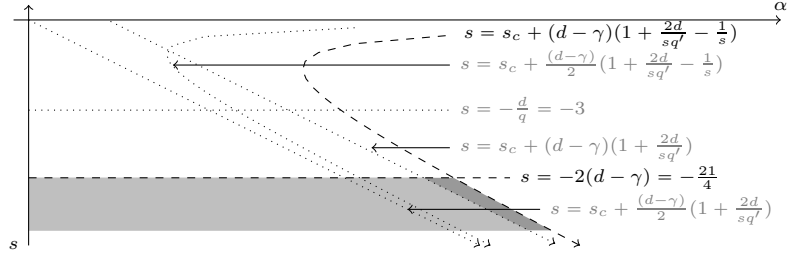
\begin{figure}
\subfigure[The case $d=3$, $\gamma=\frac{3}{8}$, $q=1$]  
{ 
\begin{tikzpicture}[scale=.4]
\fill [gray!80!white] (0,-21/4) -- plot[domain=-21/4:-7]  ({-2*\x+21/8-.25*21/\x},\x)--(15,-7)--(0,-7);
\fill [gray!50!white] (0,-21/4) -- plot[domain=-21/4:-7]  (-2*\x+21/8, {\x})--(15,-7)--(0,-7);
\draw[][->] (-.5,0)--(25,0) node[anchor=south] {\tiny{$\alpha$}};
\draw[][<-] (0,.5)--(0,-7.5) node[anchor=east] {\tiny{$s$}};
\draw[dashed][<-] [domain=-7.5:-.5] plot ({-2*\x+21/8-.25*21/\x},{\x})node[anchor=west] {\tiny{$s=s_c+(d-\gamma)(1+\frac{2d}{sq'}-\frac{1}{s})$}};
\draw[dotted][<-] [domain=17.5:21/8] plot ({\x},{-\x/2+21/16});
\draw[dotted][<-] [domain=-7.5:-.25] plot ({-2*\x-.125*21/\x},{\x});
\draw[dotted][-] [domain=0:7.5] plot ({\x},{-\x/2});
\draw[dotted][->] [domain=7.5:15] plot ({\x},{-\x/2});
\draw[dashed][-] (0,-21/4)--(16,-21/4) node[anchor= west] {\tiny{$s=-2(d-\gamma)=-\frac{21}{4}$}};
\draw[dotted][-] (0,-3)--(14,-3) node[anchor= west] {\tiny{\textcolor{gray}{$s=-\frac{d}{q}=-3$}}};
\draw[<-](4.8,-3/2)--(14,-3/2)node[anchor=west] {\tiny{\textcolor{gray}{$s=s_c+\frac{(d-\gamma)}{2}(1+\frac{2d}{sq'}-\frac{1}{s})$}}};
\draw[<-](11.4,-4.25)--(14,-4.25)node[anchor=west]{\tiny{\textcolor{gray}{$s=s_c+(d-\gamma)(1+\frac{2d}{sq'})$}}};
\draw[<-](12.6,-6.3)--(17,-6.3)node[anchor=west]{\tiny{\textcolor{gray}{$s=s_c+\frac{(d-\gamma)}{2}(1+\frac{2d}{sq'})$}}};
\end{tikzpicture}
}
\caption{Improvement by the case \ref{mt3a} in Theorem \ref{mt3} from \ref{mt11} in Theorem \ref{mt1} is shown in the darker region. }\label{f3}
\end{figure}

 \subsection{Comments on the Proof}
\begin{enumerate}
\item[--] In order to get NI based at general initial data,  we perturb general initial data by a sequence $\phi_{0,N}$ (which is defined on Fourier transform side),  to obtain $\psi_{0,N}=\psi_0+ \phi_{0,N}$.  See Subsection \ref{cp} for details.  We shall see that this $\psi_{0, N}$ will play a role of $\psi_{0,\epsilon}$ in Theorems \ref{mt1}, \ref{mt2} and \ref{mt3}.
\item[--]Kishimoto  \cite{kishimoto2019remark} proved NI   for NLS \eqref{fhe} with polynomial type nonlinearity\footnote{In \cite{kishimoto2019remark}  NLS was  treated with nonlinearity that is a polynomial in the unknown function and its complex conjugate.
}  in negative    $H^s(\rd)$ below expected critical regularity. This method  uses ideas from  Iwabuchi-Ogawa  \cite{iwabuchi2015ill} who  introduced a technique for proving ill-posedness of evolution equations, exploiting high-to-low energy transfer in the third Picard iterate, cf. \cite{bejenaru2006sharp}.  In \cite{kishimoto2019remark}, modulation space  $M^{2,1}(\rd)$  has been used as the auxiliary space  to establish norm inflation.  In \cite{iwabuchi2015ill,kishimoto2019remark},  the power series expansion of a solution to  NLS \eqref{fhe}
\[\psi_0\mapsto\psi:=\sum_{k=1}^\infty U_k[\psi_0],\]
where $U_k$ is Picard iteration  (see Definition \ref{imd}), is crucially used, see Lemma \ref{wp}  for details.  While  in \cite{oh2017remark}, Oh used  $\mathcal{F}L^1(\rd)$ as auxiliary space 
and  defined Picard iteration  $U_k$  via 
trees to establish norm inflation based at general initial data  for cubic NLS in $H^s(\mathbb R^d).$
Since $\mathcal{F}L^1(\rd)$ is a Banach algebra under pointwise multiplication,  fractional NLS \eqref{fhe}  is  locally well-posed in $\mathcal{F}L^1(\rd)$. 
On the other hand,  we do not know whether fractional Hartree equation  is wellposed  in $\mathcal{F}L^1(\R^d)$ (see \cite[Section 4]{bhimani2023hartree} for detail). However, fractional Hartree equation is globally well-posed in $\mathcal{F}L^1\cap L^2(\mathbb R^d),$ see \cite{carles2014cauchy} and \cite{BhimaniFrac}.
We shall thus use $\F L^1\cap L^2 (\R^d)$ as the auxiliary space to have the power series expansion to the solution, see Corollary \ref{uc}. 
\item[--]
The key idea to  establish norm inflation is  that the third  term in the power series expansion (in the solution given by Corollary \ref{uc})
\begin{itemize}
\item  dominates all other terms, specifically, see \eqref{dp} below and 
\item  becomes arbitrarily large at appropriate time $t=T$, see \eqref{1A} below. 
\end{itemize} 
This can be established by obtaining  the upper bound (Lemma \ref{d0}) of each component operator $U_k$ and lower bound of $U_3$ (Lemmata \ref{d2}, \ref{d2'}). 
\item[--] Due to the presence of the Coulomb type potential $$ \left| \cdot \right|^{-\gamma}$$ in the Hartree type nonlinearity, the analysis in the present paper becomes more interesting and involved.   On the Fourier side this potential becomes $\frac{c}{|\cdot|^{d-\gamma}}$.   
This potential is not in any of the Lebesgue spaces and this raises difficulties in order to get upper
 bound of Picard iterates.  
 In order to overcome this difficulty, we  divide the potential depending on the perturbation $\phi_{0,N}$.  
 More precisely,  we write $|\cdot|^{\gamma-d}=|\cdot|^{\gamma-d}\chi_{\{|\xi|\leq N^{\rho}\}}+|\cdot|^{\gamma-d}\chi_{\{|\xi|> N^{\rho}\}}$ instead of writing $|\cdot|^{\gamma-d}\chi_{\{|\xi|\leq1\}}+|\cdot|^{\gamma-d}\chi_{\{|\xi|> 1\}}$ as in the case of proving well-posedness in $\F L^1\cap L^2(\rd)$, see proof of Lemmas \ref{est0}, \ref{est}. This helps us to cover a larger range of $s$ for which NI occurs.  We then implement the refinement of techniques  from \cite{oh2017remark, kishimoto2019remark, oana,  bhimani2020norm}   to prove Theorem \ref{mt1}.  
\item[--] 
In order to prove  Theorem \ref{mt2}, we mainly follow the proof of Theorem \ref{mt1} but we crucially use the fact that at particular frequency, 
say at $n=0$ one can compare two weighted discrete Lebesgue norm: $$\n{\langle n\rangle^\sigma\F f}_{\ell^q(n=0)}=\n{\langle n\rangle^s\F f}_{\ell^q(n=0)}\quad  \text{for all} \  s,  \sigma \in \mathbb R.$$  This eventually leads to infinite loss of regularity.  See Remark \ref{psilg}  for  details.  By employing geometric  optics  and WKB analysis,  NI  with infinite loss of regularity for NLS \eqref{fhe} with  smooth gauge-invariant nonlinearities  has been established  by  
Carles and his collaborators in \cite{carles2012geometric,carles2017norm,bhimani2020norm}. 
See Appendix 2  below for details on geometric optics approach. In \cite[Proposition 1]{kishimoto2019remark}, Kishimoto proved NI with infinite loss of regularity for NLS \eqref{fhe}  with  smooth gauge-invariant nonlinearities ($|u|^{2\nu}u, \ \nu \in \mathbb N$) in $H^s$ via Fourier analytic approach. 

\item[--] 
To prove Theorem \ref{mt3}, we introduce new  resonant sets  $\mathcal{R}_{d,\alpha}$ (to be defined in Subsection \ref{res})  corresponding to fractional dispersion $(-\Delta)^{\alpha/2}$. The existence of such set allows us to improve certain estimate which we have already used  to establish Theorem \ref{mt1} c.f. Lemmas \ref{d2} and \ref{d2'}.   
Then proceeding as in the case of Theorems \ref{mt1}, \ref{mt2}, we prove Theorem \ref{mt3}. 
\end{enumerate} 
\subsection{Further Remarks}
\begin{remark} In  \cite[Theorem 1.1]{alazard2009loss}, Alazard and Carles established norm inflation with \textit{finite} loss of regularity NLS \eqref{fhe} in $H^s(\rd)$ for $0<s<s_c$  by reducing the problem to supercritical WKB analysis. We plan to address the issue of  with loss of regularity (finite and/or infinite) for \eqref{fhe} in $\mathcal{F}L^q_s(\rd)$ and $M^{p,q}_s(\rd)$ spaces when $0<s<s_c$ in a forthcoming work.
\end{remark}

\begin{remark} The analogue of  main theorems of the present paper can be proved  for  \eqref{fhe} with more general nonlinearity,  namely, $(K\ast |\psi|^{2\nu})\psi, \ \nu\in \mathbb N$, instead of $(K\ast |\psi|^2)\psi$ with straight  forward modifications in the proofs. We restrict ourselves to the $(K\ast |\psi|^2)\psi$ for the simplicity of presentation.
\end{remark}

\begin{remark} In the last decades, several authors have studied  well-posedness for wave,  Klein-Gordon  and Dirac  equations with Hartree type nonlinearity. See for example \cite{MiaoW, MiaoK, TesfahunD} and the reference therein. Our method of proof should be further applied to these and  other  nonlinear dispersive equations with Hartree type non linearity to establish norm inflation  results.
\end{remark}

\begin{remark}\label{stft}
The  STFT  of a  $f\in \mathcal{S}'(\m)$ with
 respect to a window function $0\neq g \in {\mathcal S}(\m)$ is defined by 
\begin{equation*}
V_{g}f(x,y)= \int_{\m} f(t) \overline{T_xg(t)} e^{- 2\pi i y\cdot t}dt,  \  (x, y) \in \m \times \widehat{\m}
\end{equation*}
 whenever the integral exists.  Here, $T_xg(t)=g(tx^{-1})$ is the translation operator on $\m$.  It is known \cite[Proposition 2.1]{wang2007global}, \cite{feichtinger1983modulation} that
\[ \|f\|_{M^{p,q}_s}\asymp  \left\| \|V_gf(x,y)\|_{L^p(\m)} \langle y \rangle^s \right\|_{L^q(\widehat{\m})}. \] 
 The definition of the modulation space  is independent of the choice of 
the particular window function, see  \cite[Proposition 11.3.2(c)]{GrochenigB}.
\end{remark}
\begin{remark}\label{chr} For any given function $f$ which is locally in $B$  (Banach space) (i.e,  $gf\in B, \forall g \in C_0^{\infty}(\R^d)),$ we set $f_{B}(x)= \|f g(\cdot -x)\|_{B}.$ In \cite{Fei},  Feichtinger introduced Wiener amalgam space $W(B,C)$ with local component $B$  and global component $C$ (Banach space) is defined as the space of all functions $f$ locally in $B$ such that  $f_{B}\in C$. The space $W(B, C)$ is endowed with the norm  $\|f\|_{W(B, C)}=\|f_{B}\|_{C}.$  Moreover, different choices of $g\in C^{\infty}_0(\R^d)$  generate the same space and yield equivalent norms.   We note that Fourier amalgam spaces is a Fourier image of particular Wiener amalgam spaces,  specifically,  $\mathcal{F}W(L^p, \ell^q_s)=\widehat{w}^{p,q}_s.$
\end{remark}

This article is organized as follows: In Section \ref{pr}, we gather some general tools which will be used later. In Section \ref{ki}, we establish the power series expansion of the solution along with bound on the Picard iterates in the space $\F L^1\cap L^2$. In Section \ref{th1-2}, we  prove Theorems \ref{mt1}, \ref{mt2} and in Section \ref{th3} we prove the improvement i.e. Theorem \ref{mt3}. Finally, we shall compare geometric optics and Fourier analytic approach in  Appendix 2.
 
\section{Preliminaries}\label{pr}
\noindent
\textbf{Notations}.  The notation $A \lesssim B $ means $A \leq cB$ for some constant $c > 0 $,  whereas $ A \sim B $ means $c^{-1}A\leq B\leq cA $ for some $c\geq 1$.  By $A\gtrsim B$ we mean $B\lesssim A$. For $A,B$ depending on parameter say $N\geq1$, by $A\ll B$ we mean $A/B\to0$ as $N\to\infty$ and by $A\gg B$ we mean $B\ll A$.
The Schwartz space is denoted by  $\mathcal{S}(\m)$, and the space of tempered distributions is  denoted by $\mathcal{S'}(\m).$ Let $\mathcal{F}:\mathcal{S}(\m)\to \mathcal{S}(\m)$ be the Fourier transform  defined by  
\begin{eqnarray*}
\mathcal{F}f(w)=\widehat{f}(w)=\int_{\m} f(t) e^{- 2\pi i t\cdot w}dt, \  w\in \widehat{\m}.
\end{eqnarray*}

For the convenience of reader, we recall following technical lemmata estimating a special sequence of non-negative real numbers,  convolution of two characteristic functions  and estimating the quantity $\left\|\langle \cdot \rangle^s \right\|_{L^p(Q_A)}$, $Q_A=[-A/2, A/2)^d$.

\begin{lemma}[see e.g. Lemma 4.2 in  \cite{okamoto} and Lemma 3.5 in \cite{kishimoto2019remark}]\label{it0} 
Let $\{b_k\}_{k=1}^{\infty}$ be a sequence of non-negative  real numbers such that 
\[ b_k \leq C  \sum_{{k_1,k_2, k_3 \geq
      1}\atop{k_1+k_2+ k_3 =k}}  b_{k_1}b_{k_2}
  b_{k_3} \quad\text{for all } k \geq 2.\]
Then  we have 
$ b_k \leq b_1 C_0^{k-1}$ for all $k \geq 1$, where $C_0=
  \frac{\pi^2}{2} \sqrt{C} b_1$.
\end{lemma}
\begin{lemma}\label{d3}
Let $d\geq1,  \eta_1,\eta_2\in\rd$ and $A\geq1$, then there exists $c_d,C_d>0$ such that
\begin{align*}
c_dA^d\chi_{\eta_1+\eta_2+Q_A}\leq\chi_{\eta_1+Q_A}*\chi_{\eta_2+Q_A}\leq C_dA^d\chi_{\eta_1+\eta_2+2Q_A}.
\end{align*} 
\end{lemma}
\begin{lemma}[See Lemmas 3.9 and 4.5 in \cite{bhimani2020norm}]\label{dnl} Let $A\gg 1$, $d\ge 1$, $s<0$ and $1\leq q <
  \infty$. We define 
\begin{equation*} f^q_s(A)= 
\begin{cases} 1 & \text{if} \ s< -\frac{d}{q},\\
\left( \log A \right)^{1/q} &\text{if} \ s= -\frac{d}{q},\\
A^{d/q+s} & \text{if} \  s> -\frac{d}{q}.
\end{cases}
\end{equation*}
Then we have 
$f^q_s(A) \lesssim \left\|\langle \cdot \rangle^s \right\|_{L^q(Q_A)} \sim \| (1+|n|)^s\|_{\ell^q (0\leq |n|\leq A)} \lesssim f^q_s(A)$ and  $f_s^{\infty}(A)= \left\|\langle \cdot \rangle^s \right\|_{L^{\infty}(Q_A)} \sim \left\|(1+|n|)^s \right\|_{L^{\infty}(Q_A)} \sim 1. $ 
In particular,  $f^q_s(A) \gtrsim \max( A^{\frac{d}{q}+s},1)$ for any $s<0.$
\end{lemma}
The solution to the homogeneous fractional Schr\"odinger equation with data $\psi_0$ is given by 
\begin{equation}\label{sfp}
\F U(t)\psi_0(\xi)=e^{ic_\alpha t|\xi|^{\alpha}}\F \psi_0(\xi)
\end{equation}
 where $c_\alpha=(2\pi)^\alpha, t\in \R$ and $\xi \in \widehat{\m}$.
 In the next lemma we note  quite useful fact which says the Schr\"odinger propogator $U(t)$ is uniformly  bounded on Fourier amalgam spaces.
\begin{lemma}\label{M1}
Let $t\in\R, 0<\alpha<\infty$ and $p,q\in[1,\infty]$, $s\in\R$. Then we have  
\[ \n{U(t)\psi_0}_{\widehat{w}^{p,q}_s(\m)}=\n{\psi_0}_{\widehat{w}^{p,q}_s(\m)}.
\]
\end{lemma}
\begin{proof}[{\bf Proof}] Note  that
$\n{U(t)\psi_0}_{\widehat{w}^{p,q}_s(\m)}
=\n{\n{\chi_{n+Q} e^{itc_\alpha|\xi|^\alpha}\F \psi_0}_{L^p}\langle n \rangle^s}_{\ell^q}=\n{\n{\chi_{n+Q}\F \psi_0}_{L^p}\langle n \rangle^s}_{\ell^q}.$
\end{proof}
\begin{lemma}[Algebra property]\label{m1} Let $\frac{1}{p_1}+ \frac{1}{p_2}=1+ \frac{1}{p}$ and $\frac{1}{q_1}+\frac{1}{q_2}=1+\frac{1}{q}.$ Then we have 
\[\|fg\|_{\widehat{w}^{p,q}}\lesssim \|f\|_{\widehat{w}^{p_1,q_1}}\|g\|_{\widehat{w}^{p_2,q_2}}.\]
In particular, 
$\widehat{w}^{p,q}(\m)$ is an $\F L^1$-module i.e. $\|fg\|_{\widehat{w}^{p,q}(\m)}\lesssim\|f\|_{\F L^1}\|g\|_{\widehat{w}^{p,q}(\m)}$.
\end{lemma}
\begin{proof} We note that $\mathcal{F}W(L^p, \ell^q)= \widehat{w}^{p,q}.$ Since $$L^{p_1}\ast L^{p_2}\subset L^{p} \quad \text {and} \quad \ell^{q_1}\ast \ell^{q_2}\subset \ell^{q},$$
by \cite{Fei},  we have 
$
\|fg\|_{\widehat{w}^{p,q}}=\| \hat{f}\ast \hat{g}\|_{W(L^p, \ell^q)}  \lesssim \|\hat{f}\|_{W(L^{p_1}, \ell^{q_1})}\|\hat{g}\|_{W(L^{p_2}, \ell^{q_2})} \lesssim  \|f\|_{\widehat{w}^{p_1,q_1}}\|g\|_{\widehat{w}^{p_2,q_2}}.
$
\end{proof}
\section{Key  Ingredient}\label{ki}
For $f, g, h \in \mathcal{S}(\m),$  we define \textbf{trilinear operator} associated to the nonlinearity in \eqref{fhe} as follows
\begin{equation*}
\mathcal{H}_\gamma (f,g,h)=
\begin{cases} \left[ \frac{1}{|\cdot|^\gamma}\ast\(f\bar{g}\)\right]h \quad & \text{if} \   0<\gamma <d\\
f\bar{g}h \quad  & \text{if} \ \gamma=d.
\end{cases}
\end{equation*}
When $f=g=h,$ we simply write $\mathcal{H}_{\gamma}(f,f,f)=\mathcal{H}_{\gamma}(f).$ Recall that solution to \eqref{fhe} satisfies
\begin{equation}\label{du}
\psi(t)=U(t)\psi_0-i\mu \int_0^tU(t-\tau)\mathcal{H}_\gamma(\psi(\tau))d\tau.
\end{equation}
\begin{definition}[Picard iteration]\label{imd} 
For $\psi_0\in L^2(\m),$  define
$U_1[\psi_0](t)= U(t)\psi_0$ and for $k\geq2$
\[U_k[\psi_0](t)=-i\mu\sum_{{k_1, k_2,k_3 \geq 1}\atop{k_1+k_2+k_3 =k}}\int_0^tU(t-\tau)\mathcal{H}_\gamma\left(  U_{k_1}[\psi_0],U_{k_2}[\psi_0],U_{k_3}[\psi_0] \right)(\tau) d\tau. \]
\end{definition}

\begin{remark}
The empty sums in Definition \ref{imd} are considered as zeros. In view of this one can see that $U_{2\ell}[\psi_0]=0$ for all $\ell\in\mathbb{N}$.
\end{remark}

\subsection{Local well-posedness in $L^2\cap\F L^1$}
For smooth space-time functions $\psi_1,\psi_2,\psi_3$, we set 
\begin{eqnarray}\label{nc}
\mathcal{N}(\psi_1,\psi_2,\psi_3)(t)=\int_0^tU(t-\tau)\mathcal{H}_\gamma(\psi_1(\tau),\psi_2(\tau),\psi_3(\tau))d\tau.
\end{eqnarray}
 Let $D,S$ be Banach spaces such that $S\subset C([0,T],D)$. It follows from Bejenaru and Tao's \cite[Theorem 3]{bejenaru2006sharp} that,  if the estimates
\begin{equation}\label{esti}
\n{U(\psi_0)}_{S}\lesssim\n{\psi_0}_{D},\quad\n{\mathcal{N}(\psi_1,\psi_2,\psi_3)}_{S}\lesssim\prod_{\ell=1}^3\n{\psi_\ell}_S
\end{equation}
hold, then one can find a unique solution to \eqref{du} via a Picard iteration. More precisely:

\begin{lemma}\label{wp}
Let $D,S$ be Banach spaces and $S\subset C([0,T],D)$ and $U,\mathcal{N}$ satisfy the estimates \eqref{esti}. Then there exists $\epsilon_0>0$ (depending on the constants in the inequality \eqref{esti}) such that the solution map \[\psi_0\to\psi:=\sum_{k=1}^\infty U_k[\psi_0]\] is well-defined on $B_D(\epsilon_0)=\set{\psi_0\in D:\n{\psi_0}_D\leq\epsilon_0}$ and gives a unique  solution to \eqref{du}.
\end{lemma}
\begin{lemma}\label{tle} Let $K$ be as in \eqref{hk} and $f,g, h \in \mathcal{F}L^1\cap L^2(\m).$Then we have 
$$\| (K\ast f\bar{g})h\|_{\mathcal{F}L^1\cap L^2(\m)} \lesssim \|f\|_{\mathcal{F}L^1\cap L^2 (\m)}\|g\|_{\mathcal{F}L^1 \cap L^2(\m)}\|h\|_{\mathcal{F}L^1 \cap L^2 (\m)}.$$
\end{lemma}
\begin{proof}
 By Young's inequality,  we have
\begin{eqnarray*}
\|(K\ast  (f \bar{g}))h\|_{\mathcal{F}L^p(\m)} & =& \| [\widehat{K} (\hat{f} \ast \widehat{\bar{g}})]\ast \hat{h}\|_{L^p(\widehat{\m})}\\
& \lesssim &\|\widehat{K} (\hat{f} \ast \widehat{\bar{g}})\|_{L^1(\widehat{\m})}\|h\|_{\F L^p(\m)}.
\end{eqnarray*}
Let $n>0$ and $k_1=|\cdot|^{\gamma-d}\chi_{\{|\cdot|\leq n\}}$, $k_2=|\cdot|^{\gamma-d}\chi_{\{|\cdot|>n\}}$,  where $\chi$ is the characteristic function. Then
\begin{eqnarray}
\|[|\cdot|^{\gamma-d}(\F f\ast\F\overline{g})]\|_{L^1(\widehat{\m})}
&\leq&\|k_1(\F f\ast\F\overline{g})\|_{L^1(\widehat{\m})}+\|k_2(\F f\ast\F\overline{g})\|_{L^1(\widehat{\m})}\nonumber\\
&\leq&\|k_1\|_{L^1(\widehat{\m}) }\|\F f\ast\F\overline{g}\|_{L^{\infty}(\widehat{\m})}+\|k_2\|_{L^{\infty}(\widehat{\m})}\|\F f\ast\F\overline{g}\|_{L^1(\widehat{\m})}\nonumber\\
&\lesssim& n^\gamma\|\F (f\overline{g})\|_{L^\infty(\widehat{\m})}+n^{\gamma-d}\|\F f\|_{L^1(\widehat{\m})}\|\F g\|_{L^1(\widehat{\m})}\nonumber\\
&\lesssim& n^\gamma\|f\overline{g}\|_{L^1(\m)}+n^{\gamma-d}\|\F f\|_{L^1(\m)}\|\F g\|_{L^1(\m)}\nonumber\\
&\lesssim& n^\gamma\|f\|_{L^2(\m)}\| g\|_{L^2(\m)}+n^{\gamma-d}\|f\|_{\F L^1(\m)}\|g\|_{\F L^1(\m)}.\label{estN}
\end{eqnarray}
Now the result follows from the particular case $n=1$.
\end{proof}
\begin{remark} When $\m=\mathbb R^d,$  \eqref{estN}  for $n=1$ was established in \cite{carles2014cauchy}.
Here we introduce $n$ because we shall use the estimate \eqref{estN} later with various choices of $n$'s.  This will help us to cover more range for regularity index $s$ to produce norm inflation.  We also note that  $\mathcal{F}L^1\cap L^2(\mathbb T^d)=\mathcal{F}L^1(\mathbb T^d).$  
\end{remark}
\begin{lemma}\label{cs} Let $X=\F L^1\cap L^2(\m)$ with $\|f\|_X=\max(\|f\|_{\F L^1},\|f\|_{L^2})$ and  $\psi_1, \psi_2, \psi_3 \in C([0,T], X)$.
Then
$$\|\mathcal{N}(\psi_1,\psi_2,\psi_3)\|_{C([0,T],X)}\lesssim T\prod_{\ell=1}^3\n{\psi_\ell}_{C([0,T],X)}.$$
\end{lemma}
\begin{proof} Taking \eqref{sfp}, \eqref{nc} and Lemma \ref{tle} into account, the proof follows.
\end{proof}

\subsubsection{Estimates on Picard iterates}\label{epi}
In the next lemma we establish upper bound of each Picard iterate (Definition \ref{imd}) at generic initial data $\psi_0.$
\begin{lemma}\label{est0}
Let $\psi_0\in X=\mathcal{F}L^1\cap L^2(\m)$ and $1\leq p, q \leq \infty.$ Then there exists $c>0$ such that for all $t>0$ and $k\geq3$, we have 
\begin{align}\label{wbge}
\N{U_k[\psi_0](t)}_{\widehat{w}^{p,q}(\m)}&\leq (ct)^{\frac{k-1}{2}}\N{\psi_0}_{X}^{k-1}\|\psi_0\|_{\widehat{w}^{p,q}(\m)}.
\end{align}
\end{lemma}
\begin{proof}
Let $\{b_k\}$ be a  sequence of nonnegative  real numbers such that 
\[ b_1=1\quad \text{and} \quad  b_k = \frac{4c}{k-1}  \sum_{{k_1,k_2, k_3 \geq
      1}\atop{k_1+k_2+ k_3 =k}}  b_{k_1}b_{k_2}
  b_{k_3} \quad \forall\ k \geq 2,\]
 for some $c>0$ to be chosen.  By Lemma \ref{it0}, we have $b_k\leq c_0^{k-1}$ for some $c_0>0$. Let's claim that 
\begin{align}\label{n6}
\|U_k[\psi_0](t)\|_{\widehat{w}^{p,q}(\m)}&\leq b_kt^{\frac{k-1}{2}}\|\psi_0\|_{X}^{k-1}\|\psi_0\|_{\widehat{w}^{p,q}(\m)}.
\end{align}
By Definition \ref{imd},  Lemma \ref{m1} and \eqref{estN} with $n=1$, we have  
\begin{align*}
\|U_k[\psi_0](t)\|_{\widehat{w}^{p,q}(\m)}
&\lesssim \sum_{{k_1, k_2,k_3 \geq 1}\atop{k_1+k_2+k_3 =k}}\int_0^t\big(\|U_{k_1}[\psi_0](\tau
)\|_{L^2}\|U_{k_2}[\psi_0](\tau)\|_{L^2}\\
&\hspace{2.5cm}+\|U_{k_1}[\psi_0](\tau)\|_{\F L^1}\|U_{k_2}[\psi_0](\tau)\|_{\F L^1}\big)\|U_{k_3}[\psi_0](\tau)\|_{\widehat{w}^{p,q}(\m)}d\tau.
\end{align*}
Taking $k=3$ in the above inequality,  we have 
\begin{align*}
\|U_3[\psi_0](t)\|_{\widehat{w}^{p,q}(\m)}&\leq c\int_0^t\big(\|U_1[\psi_0](\tau
)\|_{L^2}\|U_1[\psi_0](\tau)\|_{L^2}\\
&\hspace{2.5cm}+\|U_1[\psi_0](\tau)\|_{\F L^1}\|U_1[\psi_0](\tau)\|_{\F L^1}\big)\|U_1[\psi_0](\tau)\|_{\widehat{w}^{p,q}(\m)}d\tau\\
&\leq c\int_0^t\big(\|\psi_0\|_{L^2}^2+\|\psi_0\|_{\F L^1}^2\big)\|\psi_0\|_{\widehat{w}^{p,q}_s(\m)}d\tau\leq 2ct\|\psi_0\|_X^2\|\psi_0\|_{\widehat{w}^{p,q}(\m)}.
\end{align*}
Hence,   claim \eqref{n6} is true for $k=3$.  Let us assume  that \eqref{n6} is true upto $k-1$.  Then
\begin{align*}
&\|U_k[\psi_0](t)\|_{\widehat{w}^{p,q}(\m)}\\
&\leq c\sum_{{k_1, k_2,k_3 \geq 1}\atop{k_1+k_2+k_3 =k}}\int_0^tb_{k_1}b_{k_2}b_{k_3}\big(\|\psi_0\|_{L^2}^2+\|\psi_0\|_{\F L^1}^2\big)\|\psi_0\|_{X}^{k_1+k_2+k_3-3}\|\psi_0\|_{\widehat{w}^{p,q}(\m)}\tau^{\frac{k_1+k_2+k_3-3}{2}}d\tau\\
&\leq 2c\sum_{{k_1, k_2,k_3 \geq 1}\atop{k_1+k_2+k_3 =k}}b_{k_1}b_{k_2}b_{k_3}\|\psi_0\|_{X}^{k-1}\|\psi_0\|_{\widehat{w}^{p,q}(\m)}\int_0^t\tau^{\frac{k-3}{2}}d\tau\\
&\leq b_kt^{\frac{k-1}{2}}\|\psi_0\|_{X}^{k-1}\|\psi_0\|_{\widehat{w}^{p,q}(\m)}.
\end{align*}
This proves  claim \eqref{n6}.  Hence, the result follows.  This completes the proof.
\end{proof}




\begin{corollary}\label{uc}
Let $T M^2\ll1$, then  for any  $\psi_0 \in X=\F L^1\cap L^2(\m)$ with  $\|\psi_0\|_X \leq M$, there exists a unique solution $\psi\in C([0, T], X)$ to the integral equation \eqref{du} associated with \eqref{fhe}, given by 
 \begin{eqnarray}\label{pse}
 \psi= \sum_{k=1}^{\infty} U_k[\psi_0]= \sum_{\ell =0}^{\infty} U_{2 \ell +1}[\psi_0] 
 \end{eqnarray}
 which converges absolutely in $C([0, T],X).$
\end{corollary}
\begin{proof}
Taking Lemmas \ref{M1}, \ref{cs}   into account,  the proof follows as  in the proof of \cite[Lemma 2.4]{forlano2019remark} or \cite[Corollary 1]{bhimani2021norm}.
Define 
\[\Phi(\psi)(t)=U(t)\psi_0-i\mu\mathcal{N}(\psi,\psi,\psi)(t).\]
By Lemmas \ref{M1} and \ref{cs},   we have 
\begin{align*}
\n{\Phi(\psi)}_{C([0,T],X)}&\leq \n{\psi_0}_{X}+T\n{\psi}_{C([0,T],X)}^3,\\
\n{\Phi(\psi_1)-\Phi(\psi_2)}_{C([0,T],X)}&\lesssim T\max\(\n{\psi_1}_{C([0,T],X)}^2,\n{\psi_2}_{C([0,T],X)}^2\)\n{\psi_1-\psi_2}_{C([0,T],X)}.
\end{align*}
 Then considering the ball
\[B_{2M}^T=\set{\phi\in C([0,T],X):\n{\phi}_{C([0,T],X)}\leq 2M}\]
 with 
 $TM^2\ll1$ we find a fixed point of $\Phi$ in $B_{2M}^T$ and hence a solution to \eqref{du}. This  proves the existence part.
 For the second part,   in view of Lemma \ref{est0}, the series \eqref{pse} converges absolutely in $C([0,T],X)$ if $0<T\ll M^{-2}$.  In fact, for $0\le t\le T,$ using  Lemma \ref{est0} (with $p=q\in\{1,2\}$),  we  have 
\begin{eqnarray}\label{her}
\sum_{k=1}^{\infty} \|U_{k}[\psi_0](t)\|_X\le\sum_{k=1}^{\infty}(ct)^{\frac{k-1}{2}}\|\psi_0\|_X^{k}\le M\sum_{k=1}^{\infty}(\sqrt{cT}M)^{k-1}=\frac{M}{1-\sqrt{cT}M}\le2M
\end{eqnarray}
  provided $\sqrt{cT}M<1/2$.
 Then for $\epsilon>0$, there exists $j_1$ such that for all $j\geq j_1$
one has 
\begin{equation}\label{s10}
\n{\psi-\psi_j}_{C([0,T],X)}<\epsilon
\end{equation}
where \[\psi=\sum_{k=1}^\infty U_k[\psi_0] \quad\text{and}\quad \psi_j=\sum_{k=1}^j U_k[\psi_0].\]
In view of \eqref{her},  we have  $\psi,\psi_j\in B_{2M}^T$ for all $j$ as $0<T\ll M^{-2}$.    Using the continuity of $\Phi$ on $B_{2M}^T$ we find $j_2$ such that for all $j\geq j_2$ 
\begin{equation}\label{s2}
\n{\Phi(\psi)-\Phi(\psi_j)}_{C([0,T],X)}<\epsilon.
\end{equation}
Let us compute $\psi_j-\Phi(\psi_j):$
\begin{align}\label{Wr}
\psi_j-\Phi(\psi_j)&=\sum_{k=1}^j U_k[{\psi}_0]-U(t){\psi}_0+i\mu\mathcal{N}(\psi_j,\psi_j,\psi_j)\nonumber\\
&=\sum_{k=2}^j U_k[{\psi}_0]+i\mu\sum_{1\leq k_1,k_2,k_3\leq j}\mathcal{N}(U_{k_1}[{\psi}_0],U_{k_2}[{\psi}_0],U_{k_3}[{\psi}_0]).
\end{align}
Note that $U_k[{\psi}_0]=-i\mu\sum_{{k_1, k_2,k_3 \geq 1}\atop{k_1+k_2+k_3 =k}}\mathcal{N}(U_{k_1}[{\psi}_0],U_{k_2}[{\psi}_0],U_{k_3}[{\psi}_0])$ and for $1\le k\le j$ this is further equal to $-i\mu\sum_{{1\le k_1, k_2,k_3 \le j}\atop{k_1+k_2+k_3 =k}}\mathcal{N}(U_{k_1}[{\psi}_0],U_{k_2}[{\psi}_0],U_{k_3}[{\psi}_0])$. On the other hand we can rewrite
\[
\sum_{1\leq k_1,k_2,k_3\leq j}\mathcal{N}(U_{k_1}[{\psi}_0],U_{k_2}[{\psi}_0],U_{k_3}[{\psi}_0])=\sum_{k=1}^{3j} \sum_{{1\leq k_1,k_2,k_3\leq j}\atop{k_1+k_2+k_3=k}}\mathcal{N}(U_{k_1}[{\psi}_0],U_{k_2}[{\psi}_0],U_{k_3}[{\psi}_0]).
\]
Therefore,  by  \eqref{Wr} (and using the fact $\sum_{{1\le k_1, k_2,k_3 \le j}\atop{k_1+k_2+k_3 =k}}$ is a void sum for $k=1$) we  obtain
\begin{align*}
&\psi_j-\Phi(\psi_j)\\
&=-i\mu\left[\sum_{k=2}^j\sum_{{1\le k_1, k_2,k_3 \le j}\atop{k_1+k_2+k_3 =k}}\mathcal{N}(U_{k_1}[{\psi}_0],U_{k_2}[{\psi}_0],U_{k_3}[{\psi}_0])-\sum_{k=1}^{3j} \sum_{{1\leq k_1,k_2,k_3\leq j}\atop{k_1+k_2+k_3=k}}\mathcal{N}(U_{k_1}[{\psi}_0],U_{k_2}[{\psi}_0],U_{k_3}[{\psi}_0])\right]\\
&=i\mu\sum_{k=j+1}^{3j} \sum_{{1\leq k_1,k_2,k_3\leq j}\atop{k_1+k_2+k_3=k}}\mathcal{N}(U_{k_1}[{\psi}_0],U_{k_2}[{\psi}_0],U_{k_3}[{\psi}_0])=-\sum_{k=j+1}^{3j} U_{k,j}[{\psi}_0],
\end{align*}
where we set 
\[U_{k,j}[{\psi}_0]=-i\mu\sum_{{1\leq k_1,k_2,k_3\leq j}\atop{k_1+k_2+k_3=k}}\mathcal{N}(U_{k_1}[{\psi}_0],U_{k_2}[{\psi}_0],U_{k_3}[{\psi}_0]).\]
Note that $U_{k,j}$ has a less number of terms in the sum above compared to that of $U_k$.  Hence proceeding as in the proof of Lemma \ref{est}, one achieves same estimates for $U_{k,j}$.
So using $0<T\ll M^{-2},$ we obtain
\begin{align*}
\n{\psi_j-\Phi(\psi_j)}_{C([0,T],X)}&\leq\sum_{k=j+1}^{3j}\n{U_{k,j}[{\psi}_0]}_{C([0,T],X)} \leq \sum_{k=j+1}^{3j}(cT)^{\frac{k-1}{2}}\n{{\psi}_0}_{X}^k\\
&\leq M\sum_{k=j+1}^\infty(\sqrt{ cT}M)^{k-1}\leq2 M(\sqrt{cT}M)^j.
\end{align*}
Then there exists $j_3$ such that for  $j\geq j_3$ one has
\begin{equation}\label{s3}
\n{\psi_j-\Phi(\psi_j)}_{C([0,T],X)}<\epsilon.
\end{equation}
Therefore from \eqref{s10}, \eqref{s2} and \eqref{s3} one has
\[\n{\psi-\Phi(\psi)}_{C([0,T],X)}<3\epsilon.\]
Since $\epsilon>0$ is arbitrary, $\psi$ is the required fixed point for $\Phi$. 
\end{proof}
\section{Proof of Theorems \ref{mt1} and \ref{mt2}}\label{th1-2}
\subsection{Choice of perturbation to the initial base}\label{cp}
We set following notations:
\begin{itemize}
\item Let $N>0,A\geq1$ with $A\ll N,$   $Q_A=[-\frac{A}{2},\frac{A}{2})^d$ and $e_1= (1,0,...,0)\in \rd.$
\item To prove Theorems \ref{mt1} and \ref{mt2}, we  shall use
$$
\Sigma=\set{Ne_1,2Ne_1}.
$$
\item To prove  Theorem \ref{mt3},
we shall use 
$$
\Sigma=\set{Nv_1,2Nv_1,Nv_2,Nv_3},
$$
where $(v_1,v_2,v_3)\in\mathcal{R}_{d,\alpha}$ (to be defined in \eqref{rv}). The existence of such $v_1,v_2,v_3$ is guaranteed by Proposition \ref{crs}.
\item Given $\psi_0\in \widehat{w}^{p,q}_s(\m)$, we put 
\begin{eqnarray}\label{psf}
\psi_{0,N}=\psi_0+\phi_{0,N},
\end{eqnarray}
 where $\phi_{0,N}$ to be describe below (see \eqref{kid}).
 \item From \eqref{p} it follows that $\|\phi_{0,N}\|_{\w}\leq\|\phi_{0,N}\|_{\F L^1}$ 
for $N$ large enough.
\end{itemize}
Let $\psi_0\in \widehat{w}^{p,q}_s(\m)\cap\F L^1(\m)$ be given and we choose a perturbation $\phi_{0,N}$  to the initial data of the following form
 \begin{equation}\label{kid}
\F{\phi_{0,N}}(\xi)=R\chi_\Omega (\xi) \quad (\xi \in \widehat{\m})
\end{equation}
with
$$\Omega=\bigcup_{\eta\in\Sigma}(\eta+Q_A) \ \text{and} \ A=\frac{N^{\frac{a}{d}}}{(\log N)^{\frac{\tilde{a}}{d}}},$$ 
where $R=R(N),a\geq0,\tilde{a}\in\R$ to be chosen later. 
We shall see that $\psi_{0,N}$ will play a role of $\psi_{0, \epsilon}$, eventually,  in Theorems \ref{mt1}, \ref{mt2} and \ref{mt3}.
Note that  for $N$ big enough
\begin{align}\label{p}
\n{\phi_{0,N}}_{\widehat{w}^{p,q}_s(\m)}\sim RA^{d/q}N^{s}.
\end{align}
 Indeed,
\begin{align*}
 \|\phi_{0,N}\|_{\widehat{w}^{p,q}_s}=   \left\| \left\lVert \chi_{n+Q}(\xi)\F f(\xi)\right\rVert_{L_\xi^p(\mh)}  \langle n \rangle^s \right\|_{\ell^q_n(\zd)}=  R \left\| \left\lVert \chi_{n+Q}(\xi)\chi_\Omega(\xi)\right\rVert_{L_\xi^p(\mh)}  \langle n \rangle^s \right\|_{\ell^q_n(\zd)}
\end{align*}Note that $\|\chi_{n+Q}(\xi)\chi_\Omega(\xi)\|_{L_\xi^p(\mh)}$ survives only if $n\in \mathcal{G}:=\{m:(m+Q)\cap\Omega\neq\emptyset\},$ and for these $n$'s one must have $|n|\sim N$. Since $\#(\mathcal{G})\sim A^d$ and $\|\chi_{n+Q}(\xi)\chi_\Omega(\xi)\|_{L_\xi^p(\mh)}=\|\chi_{n+Q}(\xi)\|_{L_\xi^p(\mh)}=1$ for almost all $n\in\mathcal{G}$,  we conclude
\begin{align*}
 \|\phi_{0,N}\|_{\widehat{w}^{p,q}_s}&=R \left\| \left\lVert \chi_{n+Q}(\xi)\chi_\Omega(\xi)\right\rVert_{L_\xi^p(\mh)}  \langle n \rangle^s \right\|_{\ell^q_n(\zd)}= R\left(\sum_{n\in\mathcal{G}}\|\chi_{n+Q}(\xi)\chi_\Omega(\xi)\|_{L_\xi^p(\mh)}^q\langle n \rangle^{sq}\right)^{1/q}\\
 &\sim R(A^dN^{sq})^{1/q}=RA^{d/q}N^s.
\end{align*}  

\begin{remark}(density argument)  
Using the density of $Z:=\widehat{w}^{p,q}_s(\m)\cap L^2(\m)\cap\F L^1(\m)$ in $\w(\m)$, without loss of generality,  we may assume that that data $\psi_0\in Z$ 
to prove all our results in this paper. Indeed if the results are true for data  in $Z$, then for $\psi_0\in\w$ and $\varepsilon>0$,   we can choose $\psi_0^\varepsilon\in Z$ with $\|\psi_0-\psi_0^\varepsilon\|_{\w}<\varepsilon/2$. Then by hypothesis there exist $\psi_{0,\epsilon}^\varepsilon$ and $T>0$ with $\|\psi_0^\varepsilon-\psi_{0,\varepsilon}^\varepsilon\|_{\widehat{w}_{s}^{p,q}}< \varepsilon/2$, $0<T< \varepsilon$ such that the corresponding smooth solution $\psi_\varepsilon^\varepsilon$ to $\eqref{fhe}$ with data $\psi_{0,\varepsilon}^\varepsilon$ exists on $[0,T]$ and $ \|\psi_\varepsilon^\varepsilon(T)\|_{\widehat{w}^{p,q}_s}> \frac{2}{\varepsilon}>\frac{1}{\epsilon}$. Now note that $\|\psi_0-\psi_0^\varepsilon\|_{\w}+\|\psi_0^\varepsilon-\psi_{0,\varepsilon}^\varepsilon\|_{\w}<\varepsilon/2+\varepsilon/2=\varepsilon$ and set $\psi_{0,\varepsilon}:=\psi_{0,\varepsilon}^\varepsilon$, $\psi_{\varepsilon}:=\psi_{\varepsilon}^\varepsilon$.
\end{remark}
\subsection{Improved estimate for perturbation on Picard iterates}
In the next lemma we  improve  estimate \eqref{wbge} with the particular choice of data $\phi_{0,N}$.
\begin{lemma}\label{est}
Let $\phi_{0,N}$  be as in \eqref{kid}, $\rho d\leq a$ with $\tilde{a}\leq0$ or $\rho d<a$ with $\tilde{a}\in\R$. 
 Then there exists $c>0$ such that for all $t>0$ and $k\geq3$, we have 
\begin{align*}
\N{U_k[\phi_{0,N}](t)}_{\widehat{w}^{p,q}(\m)}&\leq (ctN^{\rho(\gamma-d)})^{\frac{k-1}{2}}\N{\phi_{0,N}}_{\F L^1(\m)}^{k-1}\|\phi_{0,N}\|_{\widehat{w}^{p,q}(\m)}.
\end{align*}
\end{lemma}
\begin{remark} Note that $\phi_{0,N}\in\F L^1\cap L^2(\m)$ and so Lemma \ref{est0} is applicable to it.  But Lemma \ref{est}  improves   Lemma \ref{est0} if $\rho>0$ in Hartree case in the sense that the upper bound is comparatively smaller for large $N$.
\end{remark}
\begin{proof}[Proof of Lemma \ref{est}]
Let $\{b_k\}$ be as in the proof of Lemma \ref{est0}.   We claim that 
\begin{align}\label{n7}
\|U_k[\phi_{0,N}](t)\|_{\widehat{w}^{p,q}(\m)}&\leq b_k(tN^{\rho(\gamma-d)})^{\frac{k-1}{2}}\|\phi_{0,N}\|_{\F L^1}^{k-1}\|\phi_{0,N}\|_{\widehat{w}^{p,q}(\m)}
\end{align}with appropriate $\rho$.  By Lemmata \ref{M1}, \ref{m1}, and \eqref{estN} with choosing $n=N^{\rho}$,  we have
\begin{align*}
& \|U_k[\phi_{0,N}](t)\|_{\widehat{w}^{p,q}(\m)}\\
&\lesssim \sum_{{k_1, k_2,k_3 \geq 1}\atop{k_1+k_2+k_3 =k}}\int_0^t\big(N^{\rho\gamma}\N{U_{k_1}[\phi_{0,N}](\tau)}_{L^2}\N{U_{k_2}[\phi_{0,N}](\tau)}_{L^2}\\
&\hspace{3cm}+N^{\rho(\gamma-d)}\N{ U_{k_1}[\phi_{0,N}](\tau)}_{\F L^1}\N{ U_{k_2}[\phi_{0,N}](\tau)}_{\F L^1}\big)\|U_{k_3}[\phi_{0,N}](\tau) \|_{\widehat{w}^{p,q}(\m)}d\tau.
\end{align*}
Since  $\rho d \leq a$ with $\tilde{a}\leq0$ or $\rho d<a$ with $\tilde{a}\in\R$, by \eqref{p}, we have 
$$N^{\rho\gamma}\|\phi_{0,N}\|_{L^2}^2\sim \frac{N^{\rho\gamma+a}R^2}{(\log N)^{\tilde{a}}}\leq \frac{N^{\rho(\gamma-d)+2a}R^2}{(\log N)^{2\tilde{a}}}\sim N^{\rho(\gamma-d)}\|\phi_{0,N}\|_{\F L^1}^2$$ for $N$ large.  In view of this,  we  have 
\begin{align*}
\|U_3[\phi_{0,N}](t)\|_{\widehat{w}^{p,q}(\m)}
&\leq ct\left(N^{\rho\gamma}\|\phi_{0,N}\|_{L^2}^2+N^{\rho(\gamma-d)}\|\phi_{0,N}\|_{\F L^1}^2\right)\|\phi_{0,N}\|_{\widehat{w}^{p,q}}\\
& \leq 2ctN^{\rho(\gamma-d)}\|\phi_{0,N}\|_{\F L^1}^2\|\phi_{0,N}\|_{\widehat{w}^{p,q}}.
\end{align*}
Thus claim \eqref{n7} is true for $k=3$. Let us assume the claim upto $k-1$. Then
\begin{align*}
\|U_k[\phi_{0,N}](t)\|_{\widehat{w}^{p,q}(\m)}
&\leq c\sum_{{k_1, k_2,k_3 \geq 1}\atop{k_1+k_2+k_3 =k}}\int_0^tb_{k_1}b_{k_2}(N^{\rho\gamma}\|\phi_{0,N}\|_{L^2}^{2}+N^{\rho(\gamma-d)}\|\phi_{0,N}\|_{\F L^1}^{2})\\
&\hspace{3cm}\times b_{k_3}(N^{\rho(\gamma-d)}\tau)^{\frac{k_1+k_2+k_3-3}{2}}\|\phi_{0,N}\|_{\F L^1}^{k_1+k_2+k_3-3}\N{ \phi_{0,N}}_{\widehat{w}^{p,q}(\m)}d\tau\\
&\leq b_k(N^{\rho(\gamma-d)}t)^{\frac{k-1}{2}}\|\phi_{0,N}\|_{\F L^1}^{k-1}\|\phi_{0,N}\|_{\widehat{w}^{p,q}(\m)}.
\end{align*}
Hence,  the claim is established. 
\end{proof}
 \subsection{Key Lemmas}\label{Sni}Note that Lemma \ref{d0} (2) does not work for $s>0.$ Therefore, none of our results include regularity index $s>0$.
 We shall start with following result estimating the size of the support of $\F U_k[\phi_{0,N}](t)$. This is essentially due to Kishimoto \cite{kishimoto2019remark}.


\begin{lemma}\label{kms0}
Let $a=(a_1,\cdots,a_d) \in \R^d$ and $S\subset\rd$
and $\text{dist}_\infty(a,S)=\inf_{s=(s_1,\cdots,s_d)\in S}\max_{1\le j\le d}|a_j-s_j|.$ Denote $1/2-$neighbourhood of the  set $S$ by 
$\mathbf N_{1/2}[S]=\{x\in\rd:dist_\infty(x,S)\leq1/2\}.$
There exists $C>0$ such that for  $u_0$ satisfying \eqref{kid} and $k\geq 1$, we have 
\begin{enumerate}
\item$\left|  \operatorname{supp}\F{ U_k[\phi_{0,N}]} (t)\right| \leq C^k
  A^d ,\quad \forall t\geq 0$\label{kms0a}
\item
$\left|\mathbf{N}_{1/2}  [\operatorname{supp}\F{ U_k[\phi_{0,N}]} (t)]\right| \leq 2^dC^k
  A^d ,\quad \forall t\geq 0$.\label{kms0b}
\end{enumerate}
\end{lemma}
\begin{proof}[{\bf Proof}]\eqref{kms0a}
Note that $\operatorname{supp}\F{ U_1[\phi_{0,N}]} (t)\subset\operatorname{supp}\F \phi_{0,N}$ which is contained in at most four cubes with volumes $A^d.$  Hence $\left|  \operatorname{supp}\F{ U_1[\phi_{0,N}]} (t)\right| \leq 4 A^d$ for all $t\geq 0$. It is suffice to prove that  $\operatorname{supp}\F{ U_k[\phi_{0,N}]} (t)$ is contained in $4c_d^{k-1}$ number of cubes with volume $A^d$.  Clearly the claim is true for $k=1$. Let it be true upto $k-1$ stage. Then using the property that $\operatorname{supp}  \ (f\ast g)\subset\operatorname{supp}f +\operatorname{supp} g$ one has
\begin{align}\label{6}
\operatorname{supp}\F{ U_k[\phi_{0,N}]} (t)\subset\sum_{{k_1,k_2, k_{3} \geq
      1}\atop{k_1+k_2+ k_{3} =k}} \operatorname{supp}v_{k_j}(t)
\end{align}where $v_{k_{\ell}}$ is either $\F{U_{k_{\ell}}[\phi_{0,N}]}$
or $\F{\overline{ U_{k_{\ell}}[\phi_{0,N}]}}$. Using induction we conclude that the set in RHS  of  \eqref{6} is contained in \begin{align*}
d^{2}\prod_{{k_1,k_2, k_{3} \geq
      1}\atop{k_1+k_2+ k_{3} =k}} 4c_d^{k_j-1}=4^{3}d^{2}c_d^{k-3}=4(4d)^{2}c_d^{k-3}
\end{align*}number of cubes with volume $A^d$. Set $c_d=4d$ and $C=c_d=4d$ to conclude.\\
\eqref{kms0b} 
Since $\operatorname{supp}\F{ U_k[\phi_{0,N}]} (t)$  contained in $4c_d^{k-1}$ numbers of cubes with side-lengths $A\geq1$, and  the size of 1/2-neighbourhood of such cubes is controlled by the size of the  cubes with side-lengths $2A$, the result follows.
\end{proof}

In the next lemma we establish upper bound of each Picard iterate (Definition \ref{imd}) on perturbation $\phi_{0, N}$ and on its associate $\psi_{0,N}.$
\begin{lemma}\label{d0}Let $\psi_{0, N}$  be given by \eqref{psf}, $s\leq0,  k\geq 3, \rho d\leq a$ and $1\leq p,q \leq  \infty.$ Then   \begin{enumerate}
  \item $ \|U_1[\psi_{0,N}](t)\|_{\widehat{w}^{p,q}_s(\m)} \lesssim 1+ RA^{d/q}N^s$\label{d02}
  \item $\|U_k[\phi_{0,N}](t)\|_{\widehat{w}^{p,q}_s(\m)}\leq C^k(tN^{\rho(\gamma-d)})^{\frac{k-1}{2}}(RA^d)^{k-1}R\N{\langle n\rangle^s}_{\ell^q(|n|\leq A)}$ for some $C>0.$\label{d04}
   \item \label{dpr} If  $\|\psi_0\|_Y\leq N^{\rho(\gamma-d)/2}\|\phi_{0,N}\|_Y$ with $Y=\widehat{w}^{p,q}\cap\F L^1\cap L^2$    then    
\begin{eqnarray} \label{d03}
 \|U_3[\psi_{0,N}](t)-U_3[\phi_{0,N}](t)\|_{\widehat{w}^{p,q}_s(\m)}\lesssim tR^2A^{2d}
\end{eqnarray}
 and 
 \begin{eqnarray}\label{d05}
 \|U_k[\psi_{0,N}](t)\|_{\widehat{w}^{p,q}_s(\m)}\leq C^k(tN^{\rho(\gamma-d)})^{\frac{k-1}{2}}(RA^d)^{k-1}[R\N{\langle n\rangle^s}_{\ell^q(|n|\leq A)}+\|\psi
 _0\|_YN^{\rho(d-\gamma)}]
\end{eqnarray} 
for some $C>0.$
  \end{enumerate} 
\end{lemma}
\begin{remark}  Note that hypothesis in part \eqref{dpr} is always satisfied if $\psi_0=0$. For non zero $\psi_0$, it suffices to have $N^{\rho(\gamma-d)/2}\|\phi_{0,N}\|_Y\sim RA^dN^{\rho(\gamma-d)/2+s}\gg1$ with certain choice of parameters.
\end{remark}
\begin{proof}[{\bf Proof of Lemma \ref{d0}}]  Part \eqref{d02}: By Lemma \ref{M1},  \eqref{psf} and \eqref{kid},  we have 
\begin{eqnarray*}
\|U_1[\psi_{0,N}](t)\|_{\widehat{w}^{p,q}_s(\m)} = \|\psi_{0,N}\|_{\widehat{w}^{p,q}_s(\m)} \leq \|\psi_0\|_{\w}+ RA^{d/q}N^s \lesssim 1 + RA^{d/q}N^s.
\end{eqnarray*}
Part \eqref{d04}:  For $S\subset \mathbb R^d,$ notice that 
\begin{equation}\label{N}
 \{ n\in \mathbb Z^d: (n+Q)  \cap S\neq\emptyset\} \subset \mathbf{ N}_{1/2}(S).
 \end{equation}
  In fact, $m\in  \{ n\in \mathbb Z^d: (n+Q)  \cap S\neq\emptyset\} $ implies $(m+Q)  \cap S\neq\emptyset.$   Hence, there exists $a\in S$ such that $a\in m+Q$.  As $Q=(-\frac{1}{2},\frac{1}{2}]^d$, we get $|m_j-a_j|\leq \frac{1}{2}$ for all $1\le j\le d$. Thus $dist_{\infty}(m,S)\leq\frac{1}{2}$ and so $m\in 
 \mathbf{ N}_{1/2}(S)$.   This proves \eqref{N}.   Thus we have 
\begin{align*}
\n{U_k[\phi_{0,N}](t)}_{\widehat{w}^{p,q}_s(\m)}& =\left\| \left\lVert \chi_{n+Q}(\xi)\F U_k[\phi_{0,N}](t)(\xi)\right\rVert_{L_\xi^p(\mh)}  \langle n \rangle^s \right\|_{\ell^q_n(\zd)}\\
& \leq\sup_{\xi \in \mh} |  \F{U_k[\phi_{0,N}]} (t, \xi)| \n{\<n\>^s}_{\ell^q(\{n\in\Z^d:(n+Q)\cap\ {\rm supp }\ \F U_k[\phi_{0,N}](t)\neq\emptyset\})}\\
&\leq\sup_{\xi \in \mh} |  \F{U_k[\phi_{0,N}]} (t, \xi)| \n{\<n\>^s}_{\ell^q(n\in \mathbf N_{1/2} [{\rm supp }\ \F U_k[\phi_{0,N}](t)])}\\
& \leq (ctN^{\rho(\gamma-d)})^{\frac{k-1}{2}}(RA^d)^{k-1}R\n{\<n\>^s}_{\ell^q(n\in   \mathbf{N}_{1/2} [{\rm supp }\ \F U_k[\phi_{0,N}](t)])}
\end{align*}by using Lemma \ref{est} (with $p=q=\infty$) and \eqref{p} in the final step.  
Since $s\leq 0$, for any bounded set $D\subset \R^d$, we have 
\[\| \langle n\rangle^s \|_{\ell^q (  n \in D )}\leq \| \langle n\rangle^s \|_{\ell^q ( n \in B_D )}\]
where $B_D\subset \R^d$ is the ball centered at origin with $|D|=|B_D|.$  
Using this and Lemma \ref{kms0},
 we obtain 
\begin{equation*}
\n{\<n\>^s}_{\ell^q(n\in  \mathbf{N}_{1/2}[\rm supp \ \F U_k[\phi_{0,N}](t)])}\leq
\| \langle n \rangle^s \|_{\ell^q ( \{ |n| \leq \textcolor{magenta}{2} C^{k/d} A \})}
\lesssim   C^{k/q} \|\langle n \rangle^s \|_{\ell^q(\{ |n| \leq A\})}.
\end{equation*}
Therefore 
\begin{equation}\label{s1}
\n{U_k[\phi_{0,N}](t)}_{\widehat{w}^{p,q}_s(\m)}\lesssim(ctN^{\rho(\gamma-d)})^{\frac{k-1}{2}}(RA^d)^{k-1}R\|\langle n \rangle^s \|_{\ell^q(\{ |n| \leq A\})}.
\end{equation}  This proves  part \eqref{d04}.  Next we shall prove part  \eqref{dpr}.  In view of Definition \ref{imd} and \eqref{nc}, observe that   
\begin{eqnarray*}
 I_{k}(t):&=& U_k[\psi_{0,N}](t)-U_k[{\phi}_{0,N}](t)\\
&=&-i\mu\sum_{{k_1, k_2, k_3 \geq 1}\atop{k_1+k_2+k_3=k}}\mathcal{N}(U_{k_1}[\psi_0+\phi_{0,N}],U_{k_2}[\psi_0+\phi_{0,N}],U_{k_3}[\psi_0+\phi_{0,N}])\\
&&\hspace{3cm}-\mathcal{N}(U_{k_1}[\phi_{0,N}],U_{k_2}[\phi_{0,N}],U_{k_3}[\phi_{0,N}])\\
&=&-i\mu\sum_{{k_1,k_2,k_3 \geq 1}\atop{k_1+k_2 +k_3=k}}\sum_{({\psi_1},{\psi}_{2},\psi_3)\in\mathcal{C}}\mathcal{N}(U_{k_1}[{\psi}_1],U_{k_2}[{\psi}_{2}],U_{k_2}[{\psi}_{3}])
\end{eqnarray*}
where $\mathcal{C}=\{\psi_0,\phi_{0,N}\}^3\setminus\{(\phi_{0,N},\phi_{0,N},\phi_{0,N})\}$.  
 Using the fact $s\leq 0$,  Lemma \ref{M1} and  \eqref{estN} with $n=1,$ we have
\begin{align*}
\|I_k(t)\|_{\widehat{w}_s^{p,q}} \lesssim\|I_k(t)\|_{\widehat{w}^{p,q}}%
& \lesssim \sum_{{k_1, k_2,k_3\geq 1}\atop{k_1+k_2+k_3=k}}\sum_{(v_1,v_2,v_3)\in \mathcal{C}}\int_0^t\big(\|U_{k_1}[{v}_1](\tau)\|_{\F L^1}\|U_{k_2}[v_2](\tau)\|_{{\F L^1}}\\
&\hspace{2.5cm}+\|U_{k_1}[{v}_1](\tau)\|_{ L^2}\|U_{k_2}[v_2](\tau)\|_{{L^2}}\big)\|U_{k_3}[{v}_3](\tau)\|_{\widehat{w}^{p,q}}d\tau.
\end{align*}
Using \eqref{n6} and the hypothesis,  we have 
\begin{align*}
\N{U_{k_j}[\psi_0](\tau)}_{\widehat{w}^{p,q}(\m)}&\leq b_{k_j}\tau^{\frac{k_j-1}{2}}\N{\psi_0}_{X}^{k_j-1}\|\psi_0\|_{\widehat{w}^{p,q}(\m)}\leq b_{k_j}\tau^{\frac{k_j-1}{2}}\|\psi_0\|_Y^{k_j-1}\|\psi_0\|_{\widehat{w}^{p,q}(\m)}\\
&\leq 
 b_{k_j}(\tau N^{\rho(\gamma-d)})^{\frac{k_j-1}{2}}\|\phi_{0,N}\|_Y^{k_j-1}\|\psi_0\|_{\widehat{w}^{p,q}(\m)}.
\end{align*}
On the other hand from \eqref{n7} we have 
\begin{align*}
\N{U_{k_j}[\phi_{0,N}](\tau)}_{\widehat{w}^{p,q}(\m)}\leq b_{k_j}(\tau N^{\rho(\gamma-d)})^{\frac{k_j-1}{2}}\|\phi_{0,N}\|_Y^{k_j-1}\|\phi_{0,N}\|_{\widehat{w}^{p,q}(\m)}.
\end{align*}
Since $\F L^1=\widehat{w}^{1,1}$ and $L^2=\widehat{w}^{2,2}$ we have in particular
 for each
 $(v_1,v_2,v_3)\in \mathcal{C}$,  
 \begin{eqnarray*}
 \N{U_{k_j}[v_\ell](\tau)}_{\F L^1(\m)}&\leq &b_{k_j}(\tau N^{\rho(\gamma-d)})^{\frac{k_j-1}{2}}\|\phi_{0,N}\|_Y^{k_j-1}\|v_\ell\|_{\F L^1(\m)}\\
 \N{U_{k_j}[v_\ell](\tau)}_{L^2(\m)}&\leq &b_{k_j}(\tau N^{\rho(\gamma-d)})^{\frac{k_j-1}{2}}\|\phi_{0,N}\|_Y^{k_j-1}\|v_\ell\|_{L^2(\m)}.
\end{eqnarray*}  
Combining the above inequalities,  we have 
\begin{align*}
\|I_k(t)\|_{\widehat{w}_s^{p,q}}
&\lesssim \sum_{{k_1, k_2,k_3\geq 1}\atop{k_1+k_2+k_3=k}}\sum_{(v_1,v_2,v_3)\in \mathcal{C}}\int_0^t b_{k_1}(\tau N^{\rho(\gamma-d)})^{\frac{k_1-1}{2}}\|\phi_{0,N}\|_Y^{k_1-1}\|v_1\|_Y\\
&\hspace{1cm}\times b_{k_2}(\tau N^{\rho(\gamma-d)})^{\frac{k_2-1}{2}}\|\phi_{0,N}\|_Y^{k_2-1}\|v_2\|_Y\times b_{k_3}(\tau N^{\rho(\gamma-d)})^{\frac{k_3-1}{2}}\|\phi_{0,N}\|_Y^{k_3-1}\|v_3\|_Y d\tau\\
&\lesssim N^{\rho(\gamma-d)\frac{k-3}{2}}\|\phi_{0,N}\|_Y^{k-3}\int_0^t \tau^{\frac{k-3}{2}}d\tau\sum_{{k_1, k_2,k_3\geq 1}\atop{k_1+k_2+k_3=k}}\sum_{(v_1,v_2,v_3)\in \mathcal{C}} b_{k_1}b_{k_2}b_{k_3}\|v_1\|_{Y}\|v_2\|_{Y}\|v_3\|_Y.
\end{align*}
Since each member of $\mathcal{C}$ has atleast one coordinate as $\psi_0$, we take that $\|v_j\|_{Y}$ (i.e. $\|\psi_0\|_{Y}$) out of summation, and dominate the other two terms by $\|\phi_{0,N}\|_Y$  which further leads to
\begin{align*}
\|I_k(t)\|_{\w}&\lesssim\|{\psi}_0\|_{Y} N^{\rho(\gamma-d)\frac{k-3}{2}}\|\phi_{0,N}\|_{Y}^{k-1}\int_0^t\tau^{\frac{k-3}{2}}d\tau\sum_{{k_1,k_2,k_3 \geq 1}\atop{k_1+k_2+k_3 =k}}b_{k_1} b_{k_2}b_{k_3}\\
&\leq b_kN^{\rho(\gamma-d)\frac{k-3}{2}}t^{\frac{k-1}{2}}(RA^d)^{k-1}\|{\psi}_0\|_Y\leq c^kN^{\rho(\gamma-d)\frac{k-3}{2}}t^{\frac{k-1}{2}}(RA^d)^{k-1}\|{\psi}_0\|_{Y}
\end{align*}
 as $\|\psi_0\|_Y\leq\|\phi_{0,N}\|_Y$. This proves \eqref{d03}.   By this and  \eqref{s1}  we may get  \eqref{d05}. 
\end{proof}
In the next lemma we obtain crucial lower bound for $\phi_{0,N}.$
\begin{lemma}\label{d2}  Let $\phi_{0,N}$  be given by \eqref{kid}, $1\leq p,q \leq  \infty,$ $s\in\R$, $\alpha\in(0,\infty)$ and $1\leq A \ll N$. If $0<T\ll N^{-\alpha},$ then we   have 
\[ \|U_{3} [\phi_{0,N}] (T)\|_{\widehat{w}_s^{p,q}}\geq\| \lVert \chi_{n+Q}(\xi)\F U_3[\phi_{0,N}]  (T)(\xi)\rVert_{L_\xi^p(\mh)} \<n\>^s \|_{\ell_n^q(|n|\leq A)}\gtrsim  R^3A^{2d}N^{\gamma-d}T\|\langle n \rangle^s\|_{\ell^q (|n|\leq A)}. \]
\end{lemma}
\begin{proof}[{\bf Proof}]
For notational convenience we put  $$\Gamma_\xi^+=\{(\xi_1,\xi_2):\xi_1+\xi_2=\xi\},\quad\Gamma_{\xi_1}^-=\{(\xi_3,\xi_4):\xi_3-\xi_4=\xi_1\}\quad \text{and}$$$$\Phi=c_\alpha(-|\xi|^\alpha+|\xi_2|^\alpha+|\xi_3|^\alpha-|\xi_4|^\alpha).$$
 Using  
  Definition \ref{imd} and definition of $\phi_{0,N}$ \eqref{kid},  
 we have
\begin{align}\label{dspl}
&i\mu^{-1}\F U_{3} [\phi_{0,N}] (T)(\xi) \nonumber\\
&=\int_0^Te^{ic_\alpha(T-t)|\xi|^\alpha}\left[\(\frac{1}{|\cdot|^{d-\gamma}}\(\F U_1(t)\phi_{0,N}\ast\F\overline{U_1(t)\phi_{0,N}}\)\)\ast\F U_1(t)\phi_{0,N}\right] (\xi) dt \nonumber \\
&=\int_0^Te^{ic_\alpha(T-t)|\xi|^\alpha}\left[\(\frac{1}{|\cdot|^{d-\gamma}}\(e^{ic_\alpha t|\cdot|^\alpha}\F \phi_{0,N}\ast e^{-ic_\alpha t|\cdot|^\alpha}\F \phi_{0,N}(-\cdot)\)\)\ast e^{ic_\alpha t|\cdot|^\alpha}\F\phi_{0,N}\right] (\xi) dt \nonumber\\
&=e^{ic_\alpha T|\xi|^\alpha}R^3\int_0^T\int_{\Gamma_\xi^+}\int_{\Gamma_{\xi_1}^-}\frac{e^{it\Phi}}{|\xi_1|^{d-\gamma}}\chi_{\Omega}(\xi_2)\chi_{\Omega}(\xi_3)\chi_{\Omega}(\xi_4)d\Gamma_{\xi_1}^-d\Gamma_\xi^+ dt.
\end{align}
For $\xi\in Q_A$ and $\xi_2, \xi_3, \xi_4 \in \Omega$ using $A\ll N$,  we have  $|\Phi|\lesssim N^\alpha.$  Hence, $|t\Phi|\ll 1$  for $0<t\ll N^{-\alpha}$ and so 
\[ \text{Re} \int_0^T e^{it\Phi} dt \geq \frac{T}{2}. \] 
Note that $|\xi_1|=|\xi-\xi_2|\sim N$ for $\xi \in Q_A, \xi_2\in \Omega.$
Using these observations,  Lemma \ref{d3} and \eqref{dspl},  we have 
\begin{align}\label{M2}
\abs{\F U_{3} [\phi_{0,N}] (T)(\xi)}&\gtrsim R^3T\int_{\Gamma_\xi^+}\int_{\Gamma_{\xi_1}^-}\frac{1}{|\xi_1|^{d-\gamma}}\chi_{\Omega}(\xi_2)\chi_{\Omega}(\xi_3)\chi_{\Omega}(\xi_4)d\Gamma_{\xi_1}^-d\Gamma_\xi^+ \nonumber\\
&\geq R^3TN^{\gamma-d}\int_{\Gamma_\xi^+}\int_{\Gamma_{\xi_1}^-}\chi_{\Omega}(\xi_2)\chi_{\Omega}(\xi_3)\chi_{\Omega}(\xi_4)d\Gamma_{\xi_1}^-d\Gamma_\xi^+\nonumber\\
&= R^3TN^{\gamma-d}\chi_{\Omega}\ast\chi_{(-\Omega)}\ast\chi_{\Omega}(\xi)\gtrsim R^3A^{2d}N^{\gamma-d}T\chi_{Q_A}(\xi).
\end{align}
 Now,  pointwise estimate  \eqref{M2} immediately  gives the desired estimate:
\begin{align*}
\|U_{3} [\phi_{0,N}] (T)\|_{\widehat{w}_s^{p,q}}\geq\left\| \lVert \chi_{n+Q_1}(\xi)\F U_3 [\vec{\phi}_{0,N}] (T)(\xi)\rVert_{L_\xi^p(Q_A\cap\mh)} \<n\>^{s} \right\|_{\ell_n^q} \gtrsim R^3A^{2d}N^{\gamma-d}T\|\<n\>\|_{\ell^q(|n|\leq A)}.
\end{align*}
\end{proof}

\subsection{Norm inflation  in $\widehat{w}_s^{p,q}$ }
\begin{proof}[{\bf Proof of Theorem \ref{mt1}}]
Let us first impose condition on $R,A$ so that 
\begin{equation}\label{n0}
\|\psi_0\|_{\mathcal{F}L^1\cap L^2}\leq \|\phi_{0,N}\|_{\mathcal{F}L^1\cap L^2}
\end{equation}
so that $\|\psi_{0,N}\|_X\leq\|\psi_0\|_X+\|\phi_{0,N}\|_X\sim RA^d$ (here $X=\F L^1\cap L^2(\mathcal{M})$). Then  Corollary \ref{uc} guarantees existence of the solution to  \eqref{fhe} and
the power series expansion in $X$ upto time $T$ with \begin{equation}\label{n1}
TR^2A^{2d}\ll1.
\end{equation} By Lemma \ref{d0} \eqref{dpr}, we obtain
\begin{equation}\label{series}
\N{\sum_{\ell=2}^\infty U_{2\ell+1}[\psi_{0,N}](T)}_{\widehat{w}_s^{p,q}}\leq \sum_{\ell=2}^\infty\N{U_{2\ell+1}[\psi_{0,N}](T)}_{\widehat{w}_s^{p,q}}\lesssim T^2N^{2\rho(\gamma-d)}R^5A^{4d}\|\langle \cdot \rangle^s\|_{\ell_n^q (|n|\leq A)}.
\end{equation}
provided $
TN^{\rho(\gamma-d)}R^2A^{2d}\ll1$ (which is implied by \eqref{n1} as $\rho\geq0$) and 
\begin{equation}\label{n4}
\rho d\leq{a}\text{ with }\tilde{a}\leq0\text{ or }\rho d<a\text{ with }\tilde{a}\in\R,
\end{equation}
\begin{numcases}{\|\psi_0\|_Y\leq}
 N^{\rho(\gamma-d)/2}\|\phi_{0,N}\|_Y\label{n2}\\
 N^{\rho(\gamma-d)}R\|\langle n \rangle^s\|_{\ell_n^q (|n|\leq A)}.\label{n3}
\end{numcases}
Note that \eqref{n0} follows from \eqref{n2}.  Let $\psi_N$ denote the solution of \eqref{du} corresponding to the initial data $\psi_{0,N}.$
By Corollary \ref{uc} and triangle inequality, we have 
\begin{align}
\|\psi_N(T)\|_{\w} 
&\geq  \| U_{3}[\psi_{0,N}](T)\|_{\w} - \|U_1[\psi_{0,N}](T)\|_{\w}- \left\| \sum_{\ell =2}^{\infty} U_{2\ell +1}[\psi_{0,N}](T) \right\|_{\w}.\label{14infl}
\end{align}
The main idea to  establish NI is that  the  third Picard iterate $U_3[\psi_0](T)$ dominates in power series expansion \eqref{pse}, specifically, we need
\[\|\psi_N(T)\|_{\w} \gtrsim \|U_3[\psi_{0,N}](T)\|_{\w}; \]
and so in order to ensure this, in view of \eqref{14infl},  we imposed following conditions:
\begin{equation}\label{dp}
 \|U_3[\psi_{0,N}](T)\|_{\w} \gg \|U_1[\psi_{0,N}](T)\|_{\w} \quad \text{and} \quad  \|U_3[\psi_{0,N}](t)\|_{\w} \gg  \left\| \sum_{\ell =2}^{\infty} U_{2\ell +1}[\psi_{0,N}](T) \right\|_{\w}.
\end{equation}
 We claim that it is 
enough to have the followings, as $N\to \infty,$
\begin{enumerate}
\item $RA^{\frac{d}{q}}N^s\to0$\label{i}
\item \label{ii}$TR^2A^{2d}\to0$
\item $TN^{\gamma-d}R^3A^{2d}\|\langle n \rangle^s\|_{\ell_n^q (|n|\leq A)}\to\infty$\label{iii}
\item $TN^{\gamma-d}R^3A^{2d}\|\langle n \rangle^s\|_{\ell_n^q (|n|\leq A)}\gg T^2N^{2\rho(\gamma-d)}R^5A^{4d}\|\langle n \rangle^s\|_{\ell_n^q (|n|\leq A)}$

 $\Leftrightarrow TN^{(2\rho-1)(\gamma-d)}R^2A^{2d}\ll1$\label{iv}
\item $0<T\ll N^{-\alpha}$\label{v}
\item $A\ll N$\label{vA}
\item $\|\psi_0\|_Y\ll N^{\rho(\gamma-d)/2}RA^d$\label{vi}
\item $\|\psi_0\|_Y\ll  N^{\rho(\gamma-d)}R\|\langle n \rangle^s\|_{\ell_n^q (|n|\leq A)}$\label{vii}
\item $\rho d\leq a$ with $\tilde{a}\leq0$ or $\rho d<a$ with $\tilde{a}\in\R$\label{viii} (not needed for power-type case).
\end{enumerate}
Before we proceed to establish the above conditions,   we first comment on these:
\begin{itemize}
\item  \eqref{i} ensures $\|\psi_{0,N}-\psi_0\|_{\w}\to0$ as $N\to\infty$.
\item  \eqref{ii} ensures the convergence of the infinite series in view of Corollary \ref{uc}.
\item In order to use Lemma \ref{d2} we need \eqref{v}, \eqref{vA}. 
\item In order to prove second inequality in \eqref{dp}, (using Lemma \ref{d0} \eqref{d03})  in view of Lemma \ref{d2}  and \eqref{series}, we need  \eqref{iv}. 
\item In order to prove first inequality in \eqref{dp}, we need \eqref{i} and \eqref{iii}. 
\item  \eqref{vi}, \eqref{vii} ensure \eqref{n2} and \eqref{n3}.
\item Finally \eqref{viii} ensures \eqref{n4}.
\end{itemize}
 Thus, in view of Lemma \ref{d2}, Lemma \ref{d0} \eqref{d03} and \eqref{iii},   it follows that 
\begin{equation}\label{1A} \|\psi_N(T)\|_{\w}\gtrsim \|U_3[\psi_{0,N}](T)\|_{\w} \gtrsim R^3A^{2d}N^{\gamma-d}T\|\langle n \rangle^s\|_{\ell_n^q (|n|\leq A)}\to \infty   \quad \text{as} \quad N \to \infty.\end{equation}
Hence,  the claim is established. 

{ \textbf{Cases \ref{mt11} and \ref{mt12}:} }
We shall now choose $A, R$ and $T$ as follows: 
\[
A=N^{\frac{a}{d}},\quad R=N^r\quad  \text{and} \quad T=N^{-\alpha-\epsilon}.
\]
where $a\geq0$, $r$ and $\epsilon>0$ to be chosen below.

Then using Lemma \ref{dnl}, we have $\|\langle n \rangle^s\|_{\ell^q (|n|\leq A)}\gtrsim A^{\frac{d}{q}+s}$. Therefore, it is enough to check
\begin{align*}
RA^{\frac{d}{q}}N^s&=N^{r+\frac{a}{q}+s}\ll1\\
TR^2A^{2d}&=N^{-\alpha-\epsilon+2r+2a}\ll1\\
TN^{\gamma-d}R^3A^{2d+\frac{d}{q}+s}&=N^{-\alpha-\epsilon+\gamma-d+3r+2a+\frac{a}{q}+\frac{as}{d}}\gg1\\
TN^{(2\rho-1)(\gamma-d)}R^2A^{2d}&=N^{-\alpha-\epsilon+(2\rho-1)(\gamma-d)+2r+2a}\ll1\\
TN^\alpha&=N^{-\epsilon}\ll1\\
AN^{-1}&=N^{a/d-1}\ll1\\
N^{\rho(\gamma-d)/2}RA^d&=N^{\rho(\gamma-d)/2+r+a}\gg1\text{ (not needed for }\psi_0=0)\\
N^{\rho(\gamma-d)}RA^{\frac{d}{q}+s}&=N^{\rho(\gamma-d)+r+\frac{a}{q}+\frac{as}{d}}\gg1\text{ (not needed for }\psi_0=0).
\end{align*}
Thus we only need to achieve:
\begin{enumerate}[label=(\roman*)]
\item $r+\frac{a}{q}+s<0$\label{i'}
\item $-\alpha-\epsilon+2r+2a<0$\label{ii'}
\item $-\alpha-\epsilon+\gamma-d+3r+2a+\frac{a}{q}+\frac{sa}{d}>0$\label{iii'}
\item $-\alpha-\epsilon+(2\rho-1)(\gamma-d)+2r+2a<0$\label{iv'}
\item $\epsilon>0$\label{v'}
\item $a<d$\label{vA'}
\item $\rho(\gamma-d)/2+r+a>0$ (not needed for $\psi_0=0$)\label{vi'}
\item $\rho(\gamma-d)+r+\frac{a}{q}+\frac{as}{d}>0$ (not needed for $\psi_0=0$)\label{vii'}
\item $\rho d\leq a$\label{viii'} (from \eqref{n4}, not needed for power-type case).
\end{enumerate}
Let us concentrate on the choice of $\epsilon>0$ first.  Note that above \ref{ii'},  \ref{iii'} and \ref{iv'}  conditions  are equivalent with
\begin{align*}
-\alpha+(1-2\rho)_+(d-\gamma)+2r+2a<\epsilon<-\alpha+\gamma-d+3r+2a+\frac{a}{q}+\frac{sa}{d}.
\end{align*}
To make room for $\epsilon$ we must have $
r+a\left(\frac{1}{q}+\frac{s}{d} \right)-[1+(1-2\rho)_+](d-\gamma)>0.
$ On the other hand, since we want $\epsilon>0$ (condition \ref{v'}) we also require\begin{equation}\label{A'}
-\alpha+\gamma-d+3r+2a+\frac{a}{q}+\frac{sa}{d}>0.
\end{equation}
 Thus $r$ must satisfy
\[-a\left(\frac{1}{q}+\frac{s}{d} \right)+[1+(1-2\rho)_+](d-\gamma)<r<-\frac{a}{q}-s\]where the later condition came  to above   condition  \ref{i'}.
Choose $r$ as a convex combination of its above two bounds: for $0\leq \theta\leq 1$ \begin{align}\label{r}
r&=\theta \left[-a(\frac{1}{q}+\frac{s}{d})+ 2(1-\rho)(d-\gamma)\right]-(1-\theta)\left(\frac{a}{q}+s \right)\nonumber\\
&=-a\big(\theta\frac{s}{d}+\frac{1}{q}\big)+2\theta (1-\rho) (d-\gamma)-(1-\theta)s.
\end{align}
Note that, to make room for $r$ we must have 
\[
0\leq a<d-[1+(1-2\rho)_+](d-\gamma)\frac{d}{-s}
\]which also ensures \ref{vA'}. Thus we choose
\begin{eqnarray}\label{va}
a=\delta \left( d-2(1-\rho) (d-\gamma)\frac{d}{-s}\right), \quad  (0\leq\delta<1). 
\end{eqnarray}
Since $a\geq0$, this imposes the condition\begin{equation}\label{n5}
s<-[1+(1-2\rho)_+](d-\gamma).
\end{equation} 
First let us consider $(1-2\rho)_+=1-2\rho$ i.e. $\rho\leq\frac{1}{2}$ (the case $\frac{1}{2}\leq\rho$ the corresponds to the particular case $\rho=\frac{1}{2}$). 
Putting values of  $r$ and $a$ in \eqref{A'},  we obtain
\begin{align}\label{de}
\alpha+d-\gamma
&<-3a\theta\frac{s}{d}-2\frac{a}{q}+2a+\frac{sa}{d}-3(1-\theta)s+6\theta (1-\rho) (d-\gamma)\nonumber\\
&=-a\big(3\theta\frac{s}{d}+\frac{2}{q}-2-\frac{s}{d}\big)-3(1-\theta)s+6\theta (1-\rho) (d-\gamma)\nonumber\\
&=-\delta \big( d-2(1-\rho) (d-\gamma)\frac{d}{-s}\big)\big(3\theta\frac{s}{d}+\frac{2}{q}-2-\frac{s}{d}\big)-3(1-\theta)s+6\theta (1-\rho) (d-\gamma) \nonumber\\
&=-2\delta\frac{d}{q}+2\delta d+(-s)[3\theta\delta-\delta+3(1-\theta)]+2 (1-\rho) (d-\gamma)\big[\delta\big(1-3\theta+\frac{2d}{sq'}\big)+3\theta\big]\nonumber\\
&=2\delta\frac{d}{q'}+(-s)[3\theta(\delta-1)+3-\delta]+2 (1-\rho) (d-\gamma)\big[\delta\big(1-3\theta+\frac{2d}{sq'}\big)+3\theta\big]\nonumber\\
&=:f(\theta,\delta),
\end{align} where $0<\theta<1$, $0\leq\delta<1$ are not yet fixed.
Then in order to put minimal restriction in \eqref{de}, 
we  impose the condition 
\begin{align}\label{de1}
\alpha<-(d-\gamma)+\max_{0\leq\theta,\delta\leq1}f(\theta,\delta)
\end{align}(this implies \eqref{de} with $\theta,\delta$ being maximisers (or approximate to maximizers as $\theta$ cannot take the values $0,1$ and $\delta$ cannot take the value $1$) in \eqref{de1}).

To this end,  we first rewrite 
\begin{eqnarray*}
f(\theta,\delta)=-3s-3[s+2 (1-\rho) (d-\gamma)]\left[\(\theta-\frac{1}{3}(1+\frac{2d}{sq'})\)(\delta-1)-\frac{1}{3}(1+\frac{2d}{sq'})\right].
\end{eqnarray*}
Now,  since $s+2 (1-\rho) (d-\gamma)<0$  (see \eqref{n5}), to maximize $f$ for a fixed $\rho$,  we need to maximize $\big(\theta-\frac{1}{3}(1+\frac{2d}{sq'})\big)(\delta-1)$ in side $[0,1]\times[0,1]$.  Note that in this region $\delta-1\leq0$. Therefore, $f$ will be maximum
 \begin{itemize}
 \item on $\{(\theta,\delta):\theta\in[0,1],\delta=1\}$ if $1+\frac{2d}{sq'}\leq0$ i.e. $-\frac{2d}{q'}\leq s<0$
  \item at $(\theta,\delta)=(0,0)$  if $1+\frac{2d}{sq'}\geq0$ i.e. $s\leq-\frac{2d}{q'}$.
\end{itemize}
  Therefore we will only consider the the above two values of $\theta,\delta$ ($\delta=1$ cannot be used as $\delta<1$, therefore we will use $\delta$ less but close to $1$ i.e. $0<1-\delta\ll1$).
It is clear that in terms of $\rho$, $f$ can have maximum at $\rho=0$ (for example in the first case above). On contrast \eqref{n5} is 
least restrictive for $\rho=\frac{1}{2}$ and therefore we shall consider both the cases $\rho=0$ and $\rho=\frac{1}{2}$.

We should  note that  so far we have taken care of conditions \ref{i'} to \ref{vA'}.  And we shall now verify conditions \ref{vi'} to \ref{viii'} for particular choices of $\rho.$

$\bullet$ Choose $\rho=\frac{1}{2}$ and $\delta$ so that $0<1-\delta\ll1$.  \\
With this choice of $\rho$ and $\delta,$  by continuity of $f$,  \eqref{de1} follows once we assume
\begin{align}\label{de2}
-\alpha+\gamma-d+2\frac{d}{q'}+2(-s)+(d-\gamma)(1+\frac{2d}{sq'})>0\Leftrightarrow s
&<\frac{d}{q'}-\frac{d-\gamma+\alpha}{2}+\frac{(d-\gamma)}{2}(1+\frac{2d}{sq'})\nonumber\\
&=s_c+\frac{(d-\gamma)}{2}(1+\frac{2d}{sq'}).
\end{align} or $\alpha<-2s+\frac{2d}{q'}[1+\frac{d-\gamma}{s}]$.
Conditions \ref{viii'} and \eqref{va} follows if we assume
\begin{eqnarray}\label{de3}
s<-2(d-\gamma).
\end{eqnarray}
Choosing $\theta$ near zero in \eqref{r}, from \ref{vi'} get the condition
$\frac{\gamma-d}{4} -\frac{a}{q}-s+a>0$ i.e. $\frac{\gamma-d}{4} -s+\frac{a}{q'}>0$ which is satisfied as $s<-2(d-\gamma)$.  From \ref{vii'} get  $d-\gamma\leq-2s+2(d-(d-\gamma)\frac{d}{-s})\frac{s}{d}=-2s+2s+2(d-\gamma)=2(d-\gamma)$ (which always the case).
In view of above discussion conditions \ref{i'}-\ref{viii'} satisfied if \eqref{de2} and \eqref{de3} i.e.
\begin{eqnarray}\label{de4}
s<\min \left(-2(d-\gamma),s_c+\frac{(d-\gamma)}{2}(1+\frac{2d}{sq'})\right)
\end{eqnarray}
is satisfied. 


$\bullet$ Choose $\rho=0$ and $\delta$ so that $0<1-\delta\ll1$.\\
As above  condition \eqref{de1}
 follows from $s<s_c+(d-\gamma)(1+\frac{2d}{sq'})$.  
Note that \ref{viii'} is satisfied as $\rho=0$. Choosing $\theta$ near zero, from \ref{vi'} get the condition $-\frac{a}{q}-s+a>0$ i.e. $\frac{a}{q'}-s>0$ which is always the case. From \ref{vii'} get  $0\leq -s+(d-2(d-\gamma)\frac{d}{-s})\frac{s}{d}=-s+s+(d-\gamma)=d-\gamma$ (which always the case). In view of \eqref{n5}, we have $s<-2(d-\gamma)$. Thus in view of these condition,  we have the condition
\begin{eqnarray}\label{de5}
s<\min \left(-2(d-\gamma),s_c+(d-\gamma)(1+\frac{2d}{sq'})\right).
\end{eqnarray}
Conditions \eqref{de4} and \eqref{de5} together on $s$ is nothing but the hypothesis of Theorem \ref{mt1} \ref{mt11}.  Thus,  this completes the proof of Case \ref{mt11}.\\


$\bullet$  Choose $\rho=0$, $\delta=0.$\\
In this case,  taking $\theta$ near zero:
$-\alpha+\gamma-d-3s>0\Leftrightarrow s<-\frac{d-\gamma+\alpha}{3}$. Thus the condition becomes 
\[
s<\min \left(-2(d-\gamma),-\frac{d-\gamma+\alpha}{3}\right)
\]
(\ref{vi'}, \ref{vii'} are satisfied in this case).  This completes the proof of Case \ref{mt12}.

%

{ \textbf{Case \ref{mt13}:}}

 In the proof of  cases \ref{mt11} and \ref{mt12},  we have used the estimates $f_s^q(A)\gtrsim A^{\frac{d}{q}+s}$  of Lemma \ref{dnl}.  To prove  Case \ref{mt13}     
 we instead use the estimates $f_s^q(A)\gtrsim1$.  
 We choose $A, R$ and $T$ as follows 
\[A=N^{\frac{a}{d}},\quad R=N^r,\quad T=N^{-\alpha-\epsilon}\] 
 where $\epsilon>0$. 
 Proceeding as before, to achieve NI, it is enough to satisfy the conditions \ref{i'},  \ref{ii'}, \ref{iv'}, \ref{v'}, \ref{vA'}, \ref{vi'}, \ref{viii'} along with 
 \begin{enumerate}[label=(\roman*$'$)]\setcounter{enumi}{2}
 \item\addtocounter{enumi}{4} $-\alpha-\epsilon+\gamma-d+3r+2a>0$\label{c}
 \item $\rho(\gamma-d)+r>0$ (not needed for $\psi_0=0$).\label{d}
\end{enumerate}
 Note that \ref{ii'}, \ref{iv'}, \ref{c}  are equivalent to
\begin{align}\label{a1}
-\alpha+(1-2\rho)_+(d-\gamma)+2r+2a<\epsilon<-\alpha+\gamma-d+3r+2a.
\end{align} 
To make room for $\epsilon$ in \eqref{a1},  we must have
\[
[1+(1-2\rho)_+](d-\gamma)<r.
\]
Taking $0\leq\rho\leq\frac{1}{2}$, we have $r>2(1-\rho)(d-\gamma).$ Then using this and \ref{i'}, $r$ must satisfy
\begin{equation}\label{A11}
2(1-\rho)(d-\gamma)<r<-\frac{a}{q}-s
\end{equation}and to make room for this $r$ one must have \begin{equation}\label{A10}
a<-sq-2(1-\rho)(d-\gamma)q.
\end{equation}

Assume  $s<-2(1-\rho)(d-\gamma)-\frac{d}{q}$, then $d<-sq-2(1-\rho)(d-\gamma)q$, so assume $a=\delta d$ with $\rho\leq\delta<1$ so that \ref{vA'}, \ref{viii'}, \eqref{A10} are satisfied (in power-type case, $0\leq\delta<1$). Take $r=2\theta(1-\rho)(d-\gamma)-(1-\theta)(\frac{a}{q}+s)$ with $0<\theta<1$.  In order to  have $\epsilon>0,$ from \eqref{a1} we must have
\begin{align*}
\alpha+d-\gamma&<3r+2a=-3(1-\theta)(\frac{a}{q}+s)+2a+6\theta(1-\rho)(d-\gamma)\\
&=-3s+3[\frac{d}{q}\delta+s+2(1-\rho)(d-\gamma)][\theta-1+\frac{2q}{3}]+3[s+2(1-\rho)(d-\gamma)](1-\frac{2q}{3})\\
&=:f(\theta,\delta)
\end{align*}
Now $f$ has maximum at $(\theta,\delta)=(0,\rho)$ if $q<\frac{3}{2}$ and at $(\theta,\delta)=(0,1)$ if $q\geq\frac{3}{2}$. 
Putting  $(\theta,\delta)$ close to these point of maxima, and $\rho=\frac{1}{2}$ we get the conditions
\begin{itemize}
\item $s<\frac{d}{2}(\frac{2}{3}-\frac{1}{q})-\frac{\alpha+d-\gamma}{3}$ if $q<\frac{3}{2}$ 
\item $s<{d}(\frac{2}{3}-\frac{1}{q})-\frac{\alpha+d-\gamma}{3}$ if $q\geq\frac{3}{2}$. 
\end{itemize}

{ \textbf{Case \ref{mt14}:}}
Choose $A,R$ and $T$  as follows
\[
A=\frac{N^{\frac{a}{d}}}{(\log N)^{\frac{\tilde{a}}{d}}},\quad R=N^r,\quad T=\frac{N^{-\alpha}}{(\log N)^{\epsilon}}
\]
with $\epsilon>0$ and $a\geq0, \tilde{a},r\in\R$ to be chosen below.  
We recall $f_s^q(A)= (\log A)^{1/q}$ for $s=-d/q.$
Then using \eqref{i}-\eqref{viii} it is enough to show the following
\begin{enumerate}[label={(\alph*)}]
\item\label{4a} $RA^{\frac{d}{q}}N^s=\frac{N^{\frac{1}{q}(a-d)+r}}{(\log N)^{\frac{\tilde{a}}{q}}}\to0$
\item\label{4b}  $TR^2A^{2d}=\frac{N^{-\alpha+2r+2a}}{(\log N)^{\epsilon+2\tilde{a}}}\to0$
\item\label{4c} $TN^{(1-2\rho)(d-\gamma)}R^2A^{2d}=\frac{N^{-\alpha+(1-2\rho)(d-\gamma)+2r+2a}}{(\log N)^{\epsilon+2\tilde{a}}}\to0$
\item\label{4d}   $TN^{\gamma-d}R^3A^{2d}(\log A)^{\frac{1}{q}}=\frac{N^{-\alpha+\gamma-d+2a+3r}}{(\log N)^{\epsilon+2\tilde{a}}}\left(\frac{a}{d}\log N- \frac{ \tilde{a}}{d}\log\log N\right)^{\frac{1}{q}}\to\infty$
\item\label{4e}  $TN^\alpha=\frac{1}{(\log N)^{\epsilon}}\to0$
\item\label{4eA} $AN^{-1}=\frac{N^{\frac{a}{d}-1}}{(\log N)^{\frac{\tilde{a}}{d}}}\to0$
\item\label{4f}  $N^{\rho(\gamma-d)/2}RA^d=\frac{N^{\rho(\gamma-d)/2+r+a}}{(\log N)^{\tilde{a}}}\to\infty$ (not needed for $\psi_0=0$)
\item\label{4g} $ N^{\rho(\gamma-d)}R(\log A)^{\frac{1}{q}}=N^{\rho(\gamma-d)+r}\left(\frac{a}{d}\log N- \frac{ \tilde{a}}{d}\log\log N\right)^{\frac{1}{q}}\to\infty$ (not needed for $\psi_0=0$)
\item\label{4h} $\rho d\leq a$ with $\tilde{a}\leq0$ or $\rho d<a$ with any $\tilde{a}\in\R$ (not needed for power-type case)
\end{enumerate}as $N\to\infty$.
Set $$a=\frac{1}{2}(\alpha+d-\gamma-3r).$$ Note that $0<a<d\Leftrightarrow\frac{\alpha-d-\gamma}{3}< r<\frac{\alpha+d-\gamma}{3}$ (implies \ref{4eA}). Then for \ref{4a} 
we choose \begin{align}\label{.hi}
\tilde{a}<0 \text{ with }  a+qr< d\Longleftrightarrow\alpha +(2q-3)r< d+\gamma.
\end{align} 
The condition \ref{4b}  (and hence \ref{4c} with $\rho=\frac{1}{2}$)  is satisfied
if \begin{equation}\label{.hii}
d-\gamma< r.
\end{equation} \ref{4e} is satisfied as $\epsilon>0$.
Note that 
\begin{equation}\label{.hiii}
TN^{\gamma-d}R^3A^{2d}(\log A)^{\frac{1}{q}}\gtrsim(\log N)^{\frac{1}{q}-\epsilon-2\tilde{a}}\to\infty
\end{equation}
if we choose $2\tilde{a}<\min(0,\frac{1}{q}-\epsilon)$ so that $\frac{1}{q}-\epsilon-2\tilde{a}>0$. 
This settles the condition \ref{4d}.

\fbox{Case $q<\frac{3}{2}$} One can choose $ \max\(\frac{\alpha-d-\gamma}{3},d-\gamma,\frac{d+\gamma-\alpha}{2q-3}\)< r<\frac{\alpha+d-\gamma}{3}$ to satisfy \ref{4a}, \ref{4b}, \ref{4c}. 
To make room for this $r$ we impose $2(d-\gamma)<\alpha<(\frac{3}{q}-1)d+\gamma$.
\begin{itemize}
\item Choosing $r\sim\frac{\alpha+d-\gamma}{3}$, we see \ref{4f}, \ref{4g} are satisfied. Note that this choice of $r$ will not satisfy \ref{4h}, therefore we need to choose different $r$ for the Hartree case.
\item To satisfy \ref{4h} (in Hartree case)  we impose $\alpha\leq2(q-1)\gamma-2(q-2)d$.  So that $ \max(d-\gamma,\frac{d+\gamma-\alpha}{2q-3})=d-\gamma$.   And choose $r\sim d-\gamma$ to satisfy \ref{4f}, \ref{4g}  and  impose the extra assumption $d+2(d-\gamma)\leq\alpha$ to satisfy \ref{4h}.
\end{itemize}

\fbox{Case $q\geq\frac{3}{2}$}  Choose $\max\(\frac{\alpha-d-\gamma}{3},d-\gamma\)< r<\min\(\frac{\alpha+d-\gamma}{3},\frac{d+\gamma-\alpha}{2q-3}\)$ with the assumption $2(d-\gamma)<\alpha<\min(d+\gamma,2d+2(q-1)(\gamma-d))$ (to make room for $r$).  Choose $r\sim d-\gamma$ so that \ref{4f}, \ref{4g} are satisfied. In Hartree-case the extra assumption $d+2(d-\gamma)\leq\alpha$ to satisfy \ref{4h}.\\

\textbf{
Case  \fbox{4'}}:  Note that the conditions with strict inequality follows from \ref{mt14} by putting $\gamma=d$. For the condition with equality (i.e. $\alpha=2d$) in the case $\frac{3}{2}\leq q<\infty$ we note that one can replace \eqref{.hi} by \[\tilde{a}>0 \text{ with }  a+qr\leq d\Longleftrightarrow\alpha +(2q-3)r\leq 2d\] and \eqref{.hii} by $d-\gamma\leq r$ with $\epsilon+2\tilde{a}>0$. Then by choosing $\epsilon,|\tilde{a}|>0$ small enough \eqref{.hiii} will follow if we exclude the case $q=\infty$. Now proceeding as above in \fbox{Case $q\geq\frac{3}{2}$} we achieve  the result.
 This completes the proof of part \fbox{4'} and hence Theorem \ref{mt1}.
\end{proof}

\subsection{Norm inflation with infinite loss of regularity in $\w$}   
\begin{remark}[Proof strategy for infinite loss of regularity] \label{psilg}
The idea here is to fix $A=1$ as in \cite[Proposition 1]{kishimoto2019remark}.    We shall do the same analysis as in the case of norm inflation,  so that we get the largeness of the second non-trivial Picard iteration $U_3[\psi_0]$ in $\w$.   Observing  the proof of norm inflation closely,  we  should notice that we have actually proved that the quantity $$\| \lVert \chi_{n+Q}(\xi)\F U_3[\psi_{0,N}]  (T)(\xi)\rVert_{L_\xi^p(\mh)} \<n\>^s \|_{\ell_n^q(|n|\leq A)}$$ becoming large (see \eqref{1A}).
Then we use  crucial  observation 
\begin{equation} \label{eqv}
\n{\langle n\rangle^\sigma\F f}_{\ell^q(n=0)}=\n{\langle n\rangle^s\F f}_{\ell^q(n=0)}
\end{equation}
 to establish the largeness in $\widehat{w}_\sigma^{p,q}$ for all $\sigma\in\R$.  
\end{remark}
\begin{proof}[{\bf Proof of Theorem \ref{mt2}}]
Here  we shall show that when ever we have chosen $a=\delta=\tilde{a}=0$ in the above proofs, 
we actually have infinite loss of regularity. Note that
\begin{align*}
&\|\psi_N(T)\|_{w_\sigma^{p,q}}\\
&\geq\| \lVert \chi_{n+Q}(\xi)\F \psi_N  (T)(\xi)\rVert_{L_\xi^p(\mh)} \<n\>^{\sigma} \|_{\ell_n^q(|n|\leq 1)}\\
&\sim\| \lVert \chi_{n+Q}(\xi)\F \psi_N  (T)(\xi)\rVert_{L_\xi^p(\mh)} \<n\>^s \|_{\ell_n^q(|n|\leq 1)}\\
&\geq\| \lVert \chi_{n+Q}(\xi)\F U_3[\psi_{0,N}]  (T)(\xi)\rVert_{L_\xi^p(\mh)} \<n\>^s \|_{\ell_n^q(|n|\leq 1)}- \|U_1[\psi_{0,N}](T)\|_{\w}- \left\| \sum_{\ell =2}^{\infty} U_{2\ell +1}[\psi_{0,N}](T) \right\|_{\w}
\end{align*}
Now rest of the argument goes as before starting form \eqref{14infl} as $$\| \lVert \chi_{n+Q}(\xi)\F U_3[\psi_{0,N}]  (T)(\xi)\rVert_{L_\xi^p(\mh)} \<n\>^s \|_{\ell_n^q(|n|\leq A)}$$ has same lower bound as $\|U_3[\psi_{0,N}](T)\|_{\w}$ from Lemma \ref{d2}.
\end{proof}
\section{Proof of Theorem \ref{mt3}}\label{th3}

\subsection{Improvement in presence of resonant vector} In presence of certain set of vectors referred as `resonant vector' we will improve Theorem \ref{mt1} and \ref{mt2} here. For this we define the following:
\subsection*{Resonant sets}\label{res}  We introduce resonant set $\mathcal{R}_{d,\alpha}$ corresponding to the fractional dispersion $(-\Delta)^{\alpha/2}$ as follows
\begin{eqnarray}\label{rv}
\mathcal{R}_{d,\alpha}=\set{(v_1,v_2,v_3)\in(\rd)^3:v_1,v_2\neq0,v_1+v_2=v_3,\ |v_1|^\alpha+|v_2|^\alpha=|v_3|^\alpha}.
\end{eqnarray}
The vectors in resonance set are called as resonant vector. We shall see such vector will play a crucial role to prove Theorem \ref{mt3}. We define $\mathcal{E}_d$ by  \begin{align*}
\mathcal{E}_d=\set{\alpha\in[1,\infty):\text{ there exists }(v_1,v_2,v_3)\in\mathcal{R}_{d,\alpha}}.
\end{align*}

We have imposed the condition $\alpha\geq1$ in the definition of $\mathcal{E}_d$ because for $\alpha<1$ there would not exist any $\alpha$ satisfying the defining condition of $\mathcal{E}_d$.
 This is because of the following reason: 
Note that $(a+b)^p>a^p+b^p$ for $a,b>0$ and $p>1$. Now
 $\abs{v_1}^\alpha+\abs{v_2}^\alpha=\abs{v_1+v_2}^\alpha\leq(\abs{v_1}+\abs{v_1})^\alpha$ implies
\[\abs{v_1}+\abs{v_2}\geq(\abs{v_1}^\alpha+\abs{v_2}^\alpha)^{1/\alpha}>\abs{v_1}+\abs{v_2} \text{ using }1/\alpha=p>1\]
which is a contradiction.
We note that: 
\begin{itemize}
\item $(1,1,2)\in\mathcal{R}_{1,1}$  implies that $1\in\mathcal{E}_1$,
\item for $d\geq2$, with $e_1=(1,0,\cdots,0),e_2=(0,1,0,\cdots,0)\in\rd$, we have $(e_1,e_2,e_1+e_2)\in\mathcal{R}_{d,2}$ which  implies $2\in\mathcal{E}_d$. 
\end{itemize}
In fact, our next  proposition says that  we can completely classify the set $\mathcal{E}_d.$
\begin{proposition}\label{crs} The set  $\mathcal{E}_1=\set{1}$ and for $d\geq2$, the set $\mathcal{E}_d=[1,\infty)$.
\end{proposition}
  We  postpone the proof of Proposition \ref{crs} for now and will prove in Appendix 1.  We conclude  this subsection by adding several remarks.

\begin{remark}  The resonant sets $\mathcal{R}_{d,2}$ corresponding to standard Laplacian $-\Delta$ goes back to Colliander-Keel-Staffilani-Takaoka-Tao \cite{colliander2010transfer}.   It has  played a crucial role in the previous works (via geometric optics approach) by Carles  et al.   in \cite{carles2012geometric, carles2017norm} and  by Kishimoto in (Fourier analytic approach) \cite{kishimoto2019remark}.  For details see Appendices 2 below.   Proposition \ref{crs} thus also maybe  of independent interest.
\end{remark}
Before proving the Theorem \ref{mt3}, we 
shall improve Lemma \ref{d2} by using the fact that 
there  exist non zero vectors $v_1,v_2,v_3$  in $\rd$ and satisfying the resonance condition
\[v_1+v_2=v_3\quad\text{and}\quad|v_1|^\alpha+|v_2|^\alpha=|v_3|^\alpha.\] Lemma \ref{d2'} below ensures us to have the same lower bound for $U_3[\psi_0]$ as in Lemma \ref{d2} but with lot of room to choose $T.$ The idea is to split $U_3[\psi_0]$ in two parts and have a lower bound for the first part dominating the upper bound of the second part. This splitting is inspired from 
\cite{forlano2019remark} in the context of wave equation.

\begin{lemma}\label{d2'}  Let $\psi_0$  be given by \eqref{kid},  $1\leq p,q \leq  \infty$, $s\in\R$, $\alpha\in\mathcal{E}_d$ and $1\leq A \ll N$. If $N^{-\alpha}\ll T\ll A^{-1}N^{1-\alpha},$ then we   have 
\[ \|U_{3} [\psi_0] (T)\|_{\w}\geq \lVert \chi_{n+Q}(\xi)\F U_3[\phi_{0,N}]  (T)(\xi)\rVert_{L_\xi^p(\mh)} \<n\>^s \|_{\ell_n^q(|n|\leq A)}\gtrsim  R^3A^{2d}N^{\gamma-d}T\|\langle n \rangle^s\|_{\ell^q (|n|\leq A)}
\]
\end{lemma}
\begin{proof}[{\bf Proof}]
Note that with $$\Gamma_\xi^+=\{(\xi_1,\xi_2):\xi_1+\xi_2=\xi\},\quad\Gamma_{\xi_1}^-=\{(\xi_3,\xi_4):\xi_3-\xi_4=\xi_1\},$$
$$\Phi=c_\alpha\left(-|\xi|^\alpha+|\xi_2|^\alpha+|\xi_3|^\alpha-|\xi_4|^\alpha\right)$$
we have (see \eqref{dspl})
\begin{eqnarray*}
\F U_{3} [\psi_0] (T)(\xi)
&= & e^{iT c_\alpha|\xi|^\alpha}R^3\sum_{\eta_2,\eta_3,\eta_4\in\Sigma}\int_0^T\int_{\Gamma_\xi^+}\int_{\Gamma_{\xi_1}^-}\frac{e^{it\Phi}}{|\xi_1|^{d-\gamma}}\prod_{\ell=2}^4\chi_{\eta_\ell+Q_A}(\xi_\ell)d\Gamma_{\xi_1}^-d\Gamma_\xi^+ dt\\
&= & e^{iT c_\alpha|\xi|^\alpha}R^3\sum_{(\eta_2,\eta_3,\eta_4)\in\mathcal{A}}\int_0^T\int_{\Gamma_\xi^+}\int_{\Gamma_{\xi_1}^-}\frac{e^{it\Phi}}{|\xi_1|^{d-\gamma}}\prod_{\ell=2}^4\chi_{\eta_\ell+Q_A}(\xi_\ell)d\Gamma_{\xi_1}^-d\Gamma_\xi^+ dt
\end{eqnarray*}
provided $A\ll N$ and $\xi\in Q_A$ and $\mathcal{A}=\{\eta=(\eta_2,\eta_3,\eta_4)\in\Sigma^3:\eta_2+\eta_3-\eta_4=0\}$. Set 
$$\mathcal{A}_0=\left\{\eta\in\mathcal{A}:|\eta_2|^\alpha+|\eta_3|^\alpha-|\eta_4|^\alpha=0\right\},\quad\mathcal{A}_1=\mathcal{A}\smallsetminus\mathcal{A}_0.$$
We split the sum in expression of $\F U_{3} [\psi_0] (T)(\xi)$ to write \begin{align*}
\F U_{3} [\psi_0] (T)(\xi)=I_0(\eta,T,\xi)+I_1(\eta,T,\xi)
\end{align*}
where
\begin{align*}
I_j(\eta,T,\xi)=e^{iTc_{\alpha}|\xi|^\alpha}R^3\sum_{(\eta_2,\eta_3,\eta_4)\in\mathcal{A}_j}\int_0^T\int_{\Gamma_\xi^+}\int_{\Gamma_{\xi_1}^-}\frac{e^{it\Phi}}{|\xi_1|^{d-\gamma}}\prod_{\ell=2}^4\chi_{\eta_\ell+Q_A}(\xi_\ell)d\Gamma_{\xi_1}^-d\Gamma_\xi^+ dt, \quad j=0,1.
\end{align*}
Note that $\mathcal{A}_0\subset \mathcal{R}_{d, \alpha}$ is non empty by our choice of $\Sigma$. (It is at this point we crucially use our choice of resonant sets. In fact, since $\mathcal{A}_0\neq \emptyset,$ we shall have following lower bound for $|I_0|$, which is  eventually gives the desired estimates.) 

Note that for $\eta\in\mathcal{A}_0$, $\xi_\ell\in\eta_\ell+Q_A$ one has
\begin{align*}
|\Phi|\lesssim\abs{ |\xi|^\alpha-|\eta_2|^\alpha-|\eta_3|^\alpha+|\eta_4|^\alpha}+\sum_{\ell=2}^4\abs{|\xi_\ell|^\alpha-|\eta_\ell|^\alpha}\lesssim AN^{\alpha-1}
\end{align*}

provided $\xi\in Q_A$. Therefore for $0<t<T\ll A^{-1}N^{1-\alpha}$ and $\xi\in Q_A$ we have 
\begin{align*}
\abs{I_0(\eta,T,\xi)}&\gtrsim T\int_{\Gamma_\xi^+}\int_{\Gamma_{\xi_1}^-}\frac{1}{|\xi_1|^{d-\gamma}}\prod_{\ell=2}^4\chi_{\eta_\ell+Q_A}(\xi_\ell)d\Gamma_{\xi_1}^-d\Gamma_\xi^+ \\
&\gtrsim R^3TN^{\gamma-d}\int_{\Gamma_\xi^+}\int_{\Gamma_{\xi_1}^-}\prod_{\ell=2}^4\chi_{\eta_\ell+Q_A}(\xi_\ell)d\Gamma_{\xi_1}^-d\Gamma_\xi^+ \gtrsim R^3A^{2d}N^{\gamma-d}T\chi_{Q_A}(\xi).
\end{align*}
On the other hand for $\xi\in Q_A,\eta\in\mathcal{A}_1$, $\xi_\ell\in\eta_\ell+Q_A$ we have $\abs{\Phi}\sim N^\alpha$ and $|\xi_1|=|\xi-\xi_2|\sim N$. Therefore
\begin{align*}
\abs{I_1(\eta,T,\xi)}
&=R^3\abs{\sum_{(\eta_2,\eta_3,\eta_4)\in\mathcal{A}_1}\int_{\Gamma_\xi^+}\int_{\Gamma_{\xi_1}^-}\frac{\frac{e^{it\Phi}}{\Phi}|_{t=0}^T}{|\xi_1|^{d-\gamma}}\prod_{\ell=2}^4\chi_{\eta_\ell+Q_A}(\xi_\ell)d\Gamma_{\xi_1}^-d\Gamma_\xi^+ }\\
&\leq R^3\sum_{(\eta_2,\eta_3,\eta_4)\in\mathcal{A}_1}\int_{\Gamma_\xi^+}\int_{\Gamma_{\xi_1}^-}\abs{\frac{e^{iT\Phi}-1}{\Phi}}{|\xi_1|^{\gamma-d}}\prod_{\ell=2}^4\chi_{\eta_\ell+Q_A}(\xi_\ell)d\Gamma_{\xi_1}^-d\Gamma_\xi^+ \\
&\lesssim R^3N^{\gamma-\alpha-d}\sum_{(\eta_2,\eta_3,\eta_4)\in\mathcal{A}_1}\int_{\Gamma_\xi^+}\int_{\Gamma_{\xi_1}^-}\prod_{\ell=2}^4\chi_{\eta_\ell+Q_A}(\xi_\ell)d\Gamma_{\xi_1}^-d\Gamma_\xi^+ \\
&= R^3N^{\gamma-\alpha-d}\sum_{(\eta_2,\eta_3,\eta_4)\in\mathcal{A}_1}\chi_{\eta_2+Q_A}\ast\chi_{\eta_3+Q_A}\ast\chi_{-\eta_4+Q_A}(\xi)\\
&\lesssim R^3A^{2d}N^{\gamma-\alpha-d}\chi_{Q_{3A}}(\xi),
\end{align*}
Therefore
\begin{align*}
\|U_{3} [\psi_0] (T)\|_{\mathcal{F}L_s^p}\geq\| \lVert \chi_{n+Q_1}(\xi)\F U_3[\phi_{0,N}]  (T)(\xi)\rVert_{L_\xi^p(\mh)} \<n\>^s \|_{\ell_n^q(|n|\leq A)}\gtrsim  R^3A^{2d}N^{\gamma-d}T\|\langle n \rangle^s\|_{\ell^q (|n|\leq A)}
\end{align*}
provided $N^{-\alpha}\ll T\ll A^{-1}N^{1-\alpha}$.
\end{proof}

\begin{proof}[{\bf Proof of Theorem \ref{mt3}}]
Taking  Lemmata \ref{d2'} and \ref{d2}  into account, we note that we can replace the condition \eqref{v} in the proof of Theorem \ref{mt1} by \begin{equation}\label{5'}
\quad 0< T\ll A^{-1}N^{1-\alpha}\quad\text{ and }\quad T\not\sim N^{-\alpha}
\end{equation} 
{ \textbf{Cases \ref{mt3a} and \ref{mt3b}:}}
Taking 
\[A=N^{\frac{a}{d}},\quad R
=N^r,\quad T=N^{-\frac{a}{d}+1-\alpha-\epsilon}\]
with $0<\epsilon\neq1$. 
Note that the only difference compared to the proofs of previous results lies in terms of power of $T$. Therefore we get all the results as before with minor changes. We have to check the conditions that involves $T$ i.e. conditions \eqref{ii}, \eqref{iii}, \eqref{iv}, \eqref{v}. 
For this, precisely we need to replace $-\alpha$  by $-\frac{a}{d}+1-\alpha$ in the power of $N$. 
Therefore proceeding as in the proof of Theorem \ref{mt1} we require
\begin{equation}
s<-[1+(1-2\rho)_+](d-\gamma)
\end{equation} 
and 
\begin{align*}
-1+\alpha+d-\gamma&<-\frac{a}{d}-a\big(3\theta\frac{s}{d}+\frac{2}{q}-2-\frac{s}{d}\big)-3(1-\theta)s+6\theta (1-\rho) (d-\gamma)\\
&=-a\big(3\theta\frac{s}{d}+\frac{2}{q}-2+\frac{1}{d}-\frac{s}{d}\big)-3(1-\theta)s+6\theta (1-\rho) (d-\gamma)\\
&=-3s-3[s+2 (1-\rho) (d-\gamma)]\\
&\quad\times\left[\(\theta-\frac{1}{3}(1+\frac{2d}{sq'}-\frac{1}{s})\)(\delta-1)-\frac{1}{3}(1+\frac{2d}{sq'}-\frac{1}{s})\right]
\end{align*}
As before we shall put various values of $\theta,\delta,\rho$ to have various conditions.

$\bullet$ Choose $\rho=\frac{1}{2}$ and $\delta$ so that $0<1-\delta\ll1$.  So that the above condition follows from
\begin{align*}
&1-\alpha+\gamma-d-3s+[s+ (d-\gamma)](1+\frac{2d}{sq'}-\frac{1}{s})>0\\
\Longleftrightarrow&-\alpha+\gamma-d-2s+\frac{2d}{sq'}+(d-\gamma)](1+\frac{2d}{sq'}-\frac{1}{s})>0\\
\Longleftrightarrow& s<\frac{d}{q'}-\frac{d-\gamma+\alpha}{2}+\frac{(d-\gamma)}{2}(1+\frac{2d}{sq'}-\frac{1}{s})=s_c+\frac{(d-\gamma)}{2}(1+\frac{2d}{sq'}-\frac{1}{s}).
\end{align*}Thus the condition becomes 
\[
s<\min\(-2(d-\gamma),s_c+\frac{(d-\gamma)}{2}(1+\frac{2d}{sq'}-\frac{1}{s})\).
\]

$\bullet$ Choose $\rho=0$ and $\delta$ so that $0<1-\delta\ll1$. As before we get
\[
s<\min\(-2(d-\gamma),s_c+(d-\gamma)(1+\frac{2d}{sq'}-\frac{1}{s})\).
\]

$\bullet$  Choosing $\rho=0$, $\delta=0$ and $\theta$ near zero:
$1-\alpha+\gamma-d-3s>0\Leftrightarrow s<-\frac{d-\gamma+\alpha-1}{3}$. Thus the condition becomes 
\[
s<\min \left(-2(d-\gamma),-\frac{d-\gamma+\alpha-1}{3}\right)
\]


$\bullet$ Choosing $\rho,\delta=0$, $\theta$ near $1$:
$1-\alpha+\gamma-d+6(d-\gamma)>0\Leftrightarrow \frac{\alpha-1}{5}+\gamma<d$. Thus the condition becomes 
\[
s<-2(d-\gamma),\quad\frac{\alpha-1}{5}+\gamma<d.
\]

{\textbf{Case \ref{mt3c}:}}
Proceeding as in the  proof of  case \ref{mt13} 
of Theorem \ref{mt1} and using Lemma \ref{d2'} as above, $-\frac{a}{d}+1-\alpha+\gamma-d+3r+2a>0$ we get the condition $s<-(d-\gamma)-\frac{d}{q}$ with \begin{itemize}
\item $s<\frac{d}{2}(\frac{2}{3}-\frac{1}{3d}-\frac{1}{q})-\frac{\alpha-1+d-\gamma}{3}$ if $q\leq\frac{3d}{2d-1}$
\item $s<{d}(\frac{2}{3}-\frac{1}{3d}-\frac{1}{q})-\frac{\alpha-1+d-\gamma}{3}$ if $q\geq\frac{3d}{2d-1}$ (this does not give new result c.f. case \ref{mt13} 
of Theorem \ref{mt1}).
\end{itemize}

{\textbf{Case \ref{mt3d}:}}
Choose $A,R,T$ in the following way
\[
A=\frac{N^{\frac{a}{d}}}{(\log N)^{\frac{\tilde{a}}{d}}},\quad R=N^r,\quad T=\frac{N^{-\frac{a}{d}+1-\alpha}}{(\log N)^{\epsilon-\frac{\tilde{a}}{d}}}
\]
with $\epsilon>0$ and $a, \tilde{a}$ to be chosen below. We recall $f_s^q(A)= (\log A)^{1/q}$ for $s=-d/q.$
Then as in part \ref{mt14} we need to have
\begin{enumerate}[label={(\alph*')}]
\item\label{.a} $RA^{\frac{d}{q}}N^s=\frac{N^{\frac{1}{q}(a-d)+r}}{(\log N)^{\frac{\tilde{a}}{q}}}\to0$
\item\label{.b}  $TR^2A^{2d}=\frac{N^{-\frac{a}{d}+1-\alpha+2a+2r}}{(\log N)^{\epsilon-\frac{\tilde{a}}{d}+2\tilde{a}}}\to0$
\item\label{.c}   $TN^{(1-2\rho)(d-\gamma)}R^2A^{2d}=\frac{N^{-\frac{a}{d}+1-\alpha+(1-2\rho)(d-\gamma)+2r}}{(\log N)^{\epsilon-\frac{\tilde{a}}{d}+2\tilde{a}}}\to0$
\item\label{.d}   $TN^{\gamma-d}R^3A^{2d}(\log A)^{\frac{1}{q}}=\frac{N^{-\frac{a}{d}+1-\alpha+\gamma-d+2a+3r}}{(\log N)^{\epsilon-\frac{\tilde{a}}{d}+2\tilde{a}}}\left(\frac{a}{d}\log N- \frac{ \tilde{a}}{d}\log\log N\right)^{\frac{1}{q}}\to\infty$
\item\label{.e}   $TAN^{\alpha-1}=\frac{1}{(\log N)^{\epsilon}}\to0$
\item\label{.f} $N^{\rho(\gamma-d)/2}RA^d=\frac{N^{\rho(\gamma-d)/2+r+a}}{(\log N)^{\tilde{a}}}\to\infty$ (not needed for $\psi_0=0$)
\item\label{.g} $ N^{\rho(\gamma-d)}R(\log A)^{\frac{1}{q}}=N^{\rho(\gamma-d)+r}\left(\frac{a}{d}\log N- \frac{ \tilde{a}}{d}\log\log N\right)^{\frac{1}{q}}\to\infty$ (not needed for $\psi_0=0$)
\item\label{.h} $\rho d\leq a$ with $\tilde{a}\leq0$ (not needed for power-type case)
\end{enumerate}
as $N\to\infty$.
Set $a=(2-\frac{1}{d})^{-1}(\alpha-1+d-\gamma-3r)$, note that $0<a<d \Leftrightarrow\frac{\alpha-d-\gamma}{3} <r<\frac{\alpha-1+d-\gamma}{3}$. Then for \ref{.a} we impose 
$\tilde{a}<0$ and \[a+qr< d\Longleftrightarrow\alpha +((2-\frac{1}{d})q-3)r< d+\gamma.\]
For \ref{.b} (and hence \ref{.c} with $\rho=\frac{1}{2}$)  we impose 
$d-\gamma< r$. 
For \ref{.d}
\[
TN^{\gamma-d}R^3A^{2d}(\log A)^{\frac{1}{q}}\gtrsim(\log N)^{\frac{1}{q}-\epsilon+\frac{\tilde{a}}{d}-2\tilde{a}}\to\infty
\]
if we choose $\tilde{a}$ small enough so that $\frac{1}{q}-\epsilon+\frac{\tilde{a}}{d}-2\tilde{a}>0$. \ref{.e} holds as $\epsilon>0$.

\fbox{ if $q\leq\frac{3d}{2d-1}$} one can choose $ \max(\frac{\alpha-d-\gamma}{3},d-\gamma,\frac{d+\gamma-\alpha}{(2-1/d)q-3})< r<\frac{\alpha-1+d-\gamma}{3}$ to satisfy the above three conditions. To make room for this $r$ we impose $2(d-\gamma)+1<\alpha< (\frac{6d}{(2d-1)q}-1)d+\gamma+(1-\frac{3d}{(2d-1)q})$. Note that \ref{.f}, \ref{.g} are satisfied for $r=d-\gamma$ and hence for $r=\frac{\alpha-1+d-\gamma}{3}>d-\gamma$.

\begin{itemize}
\item choosing $r\sim\frac{\alpha-1+d-\gamma}{3}$ \ref{.f}, \ref{.g} are satisfied. Note that this choice does satisfy \ref{.h}.

\item Assume $\alpha<[4-(2-\frac{1}{d})q]d-[2-(2-\frac{1}{d})q]\gamma$ so that we can choose $r\sim d-\gamma$ (for which \ref{.f}, \ref{.g} are satisfied). For Hartree case we farther need $\alpha>3d-2\gamma+1$ to satisfy \ref{.h}.
\end{itemize}

\fbox{if $q>\frac{3d}{2d-1}$} choose $\max(\frac{\alpha-d-\gamma}{3},d-\gamma)<r<\min(\frac{\alpha-1+d-\gamma}{3},\frac{d+\gamma-\alpha}{(2-1/d)q-3})$ with the assumption $1+2(d-\gamma)<\alpha<\min(d+\gamma,4-(2-\frac{1}{d})q]d-[2-(2-\frac{1}{d})q]\gamma)$.    Choosing $r\sim d-\gamma$, \ref{.f} and \ref{.g} are satisfied.   For Hartree case we farther need $\alpha>3d-2\gamma+1$ to satisfy \ref{.h}. 

The case $\alpha=2d$ in \fbox{d'} is similarly follow by modifying argument for \ref{mt3d} as in the proof \ref{mt14} of Theorem \ref{mt1}.
\end{proof}

\section*{Appendix }
\noindent
\subsection*{1. Characterization of resonant sets}.
In this section we shall classify  resonant sets $\mathcal{R}_{d,\alpha}$ and find the set $\mathcal{E}_d$, defined in Subsection \ref{res}.
\begin{L1}\label{A0}
The set $\mathcal{R}_{1,1}=\set{(aN,(1-a)N,N):0\leq a\leq1, N\in\R}$ and for $\alpha\in(0,\infty)\setminus\{1\}$
\[\mathcal{R}_{1,\alpha}=\set{(N,0,N),(0,N,N):N\in\R}.\]
\end{L1}
\begin{proof}[{\bf Proof}]
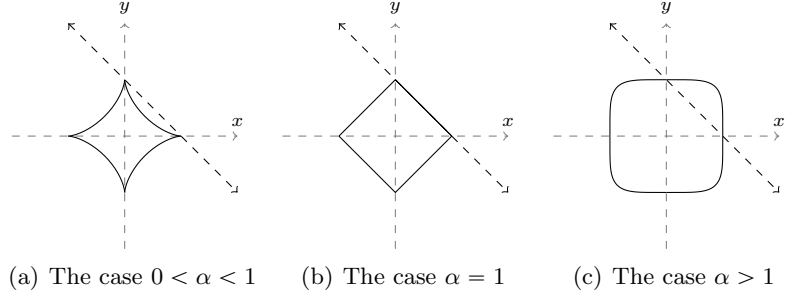
\begin{figure}
\subfigure[The case $0<\alpha<1$]  
{ 
\begin{tikzpicture}[scale=.75]
\draw[dashed,gray][->] (-2,0)--(2,0);
\draw[dashed,gray][->] (0,-2)--(0,2);
\draw[dashed][<->] (-1,2)--(2,-1);
\draw [domain=0:90] plot ({(cos(\x))^(3)}, {(sin(\x))^(3)});
\draw [domain=90:180] plot ({-(-cos(\x))^(3)}, {(sin(\x))^(3)});
\draw [domain=180:270] plot ({-(-cos(\x))^(3)}, {-(-sin(\x))^(3)});
\draw [domain=270:360] plot ({(cos(\x))^(3)}, {-(-sin(\x))^(3)});
\draw[](2,0)node[anchor=south]{\tiny{$x$}};
\draw[](0,2)node[anchor=south]{\tiny{$y$}};
\end{tikzpicture}
}
\subfigure[The case $\alpha=1$]  
{ 
\begin{tikzpicture}[scale=.75]
\draw[dashed,gray][->] (-2,0)--(2,0); 
\draw[dashed,gray][->] (0,-2)--(0,2);
\draw[dashed][<->] (-1,2)--(2,-1);
\draw[][-] (1,0)--(0,1);
\draw[][-] (1,0)--(0,1)--(-1,0)--(0,-1)--(1,0);
\draw[](2,0)node[anchor=south]{\tiny{$x$}};
\draw[](0,2)node[anchor=south]{\tiny{$y$}};
\end{tikzpicture}
}
\subfigure[The case $\alpha>1$]  
{ 
\begin{tikzpicture}[scale=.75]
\draw[dashed,gray][->] (-2,0)--(2,0);
\draw[dashed,gray][->] (0,-2)--(0,2);
\draw[dashed][<->] (-1,2)--(2,-1);
\draw [domain=0:90] plot ({(cos(\x))^(2/5)}, {(sin(\x))^(2/5)});
\draw [domain=90:180] plot ({-(-cos(\x))^(2/5)}, {(sin(\x))^(2/5)});
\draw [domain=180:270] plot ({-(-cos(\x))^(2/5)}, {-(-sin(\x))^(2/5)});
\draw [domain=270:360] plot ({(cos(\x))^(2/5)}, {-(-sin(\x))^(2/5)});
\draw[](2,0)node[anchor=south]{\tiny{$x$}};
\draw[](0,2)node[anchor=south]{\tiny{$y$}};
\end{tikzpicture}
}
\caption{\small{Graph of $x+y=v_3$ (black dashed line) and $|x|^\alpha+|y|^\alpha=|v_3|^\alpha$ (black line) with various $\alpha>0$ in the case $v_3>0$.
 }
 }\label{f5}
\end{figure}

Although the proof is clear from Figure \ref{f5}, here we shall present an analytic proof.
Let $(v_1,v_2,v_3)\in\mathcal{R}_{1,\alpha}$. Then $(v_1,v_2)$ has to solve the system of equations
\begin{align*} 
\begin{rcases}
\hspace{.9cm}x+y=v_3\\
|x|^\alpha+|y|^\alpha=|v_3|^\alpha.
\end{rcases}
\end{align*} 
Note that $x+y=v_3$ has the solutions $(av_3,(1-a)v_3,v_3)$ with $a\in\R$. Now $|x|^\alpha+|y|^\alpha=|v_3|^\alpha$ implies $a$ must satisfy 
\begin{equation*}\label{Ap1}
\abs{a}^\alpha+\abs{1-a}^\alpha=1.
\end{equation*}  
This enforces $\abs{a}\leq1,\abs{1-a}\leq1$ which implies $0\leq a\leq1$.
In the case $\alpha=1$, the above equality is satisfied by all $0\leq a\leq1$. But for $\alpha\in(0,\infty)\setminus\{1\}$, the above equality is satisfied only by $a=0,1$.
Indeed if $\alpha>1$ and $0<a<1$, then \[\abs{a}^\alpha+\abs{1-a}^\alpha<a +1-a=1\] andif $\alpha<1$ and $0<a<1$, then  $1=a+1-a<\abs{a}^\alpha+\abs{1-a}^\alpha$ which are contradictions.
\end{proof}

\begin{L2}\label{A1} 
Let $\alpha\geq1$, $d\geq2$ and set $v_1=(\cos^{2/\alpha} \theta,0,0,\cdots,0)$,
\begin{align*}
v_2=\bigg(\frac{1+\cos^{4/\alpha}\theta-\sin^{4/\alpha}\theta}{2\cos^{2/\alpha}\theta}-\cos^{2/\alpha} \theta,\frac{\sqrt{2(\cos^{4/\alpha}\theta+\sin^{4/\alpha}\theta)-(\cos^{4/\alpha}\theta-\sin^{4/\alpha}\theta)^2-1}}{\pm2\cos^{2/\alpha}\theta},\\
0,\cdots,0\bigg)
\end{align*}
with $\theta\in[0,\frac{\pi}{2})$. Then,
\begin{enumerate}
\item $(v_1,v_2,v_1+v_2)\in\mathcal{R}_{d,\alpha}$.\label{A1a}
\item \label{A1b} Upto redefinition of coordinate, any non-zero element of $\mathcal{R}_{d,\alpha}$ 
is of the form $N(v_1,v_2,v_3)$ where $N\in(0,\infty)$.
\end{enumerate}
\end{L2}
\begin{proof}[{\bf Proof}]
\eqref{A1b}
Let $(w_1,w_2,w_1+w_2=w_3)\in\mathcal{R}_{d,\alpha}$. Without loss of generality assume $|w_3|=1$ (otherwise one can normalise it by multiplying by $|w_3|^{-1}$), then after redefining the coordinate if needed, 
one can assume that $w_1,w_2,w_3$ are of the form \[w_1=(a_1,0,0,\cdots,0),\quad w_2=(a_2,b_2,0,\cdots,0)\quad\text{and}\quad w_3=(\cos r,\sin r,0,\cdots,0)\] with $0\leq r<2\pi$.
Then we have $a_1+a_2=\cos r, b_2=\sin r$ and 
\[|w_2|^2=a_2^2+b_2^2=(a_1-\cos r)^2+\sin^2r=1-2a_1\cos r+a_1^2.\]
Since $|w_1|^\alpha+|w_2|^\alpha=|w_3|^\alpha=1$, one has
\[|a_1|^\alpha+(1-2a_1\cos r+a_1^2)^{\alpha/2}=1\]
Then $a_1, r$ has to satisfy $|a_1|^{\alpha/2}=\cos\theta$ and $(1-2a_1\cos r+a_1^2)^{\alpha/4}=\sin\theta$ for some $\theta\in[0,\frac{\pi}{2})$. This implies \[a_1=(\cos\theta)^{2/\alpha}\quad \text{and}\quad 1-2a_1\cos r+a_1^2=(\sin\theta)^{4/\alpha}.\] Thus $\cos r$ is given by
\[\cos r=\frac{1+a_1^2-\sin^{4/\alpha}\theta}{2a_1}=\frac{1+\cos^{4/\alpha}\theta-\sin^{4/\alpha}\theta}{2\cos^{2/\alpha}\theta}.\]
Note that, to make sense the above inequality, we must have the RHS in the interval $[-1,1]$. 
This is equivalent to
\begin{align*}
1+\cos^{4/\alpha}\theta-\sin^{4/\alpha}\theta\leq 2\cos^{2/\alpha}\theta\Longleftrightarrow
(1-\cos^{2/\alpha})^2\leq\sin^{4/\alpha}\theta
\Longleftrightarrow1\leq \cos^{2/\alpha}\theta+\sin^{2/\alpha}\theta\Longleftrightarrow\alpha\geq1
\end{align*} and
\[-2\cos^{2/\alpha}\theta\leq1+\cos^{4/\alpha}\theta-\sin^{4/\alpha}\Longleftrightarrow\sin^{2/\alpha}\theta\leq
1+\cos^{2/\alpha}\theta.
\]
Note that the last condition is satisfied as $\cos^{2/\alpha}\theta\geq0$.
Then \begin{align*}
\sin^2r&=1-\frac{1+\cos^{8/\alpha}\theta+\sin^{8/\alpha}\theta+2\cos^{4/\alpha}\theta-2\sin^{4/\alpha}\theta-2\cos^{4/\alpha}\theta\sin^{4/\alpha}\theta}{4\cos^{4/\alpha}\theta}\\
&=\frac{2\cos^{4/\alpha}\theta+2\cos^{4/\alpha}\theta\sin^{4/\alpha}\theta+2\sin^{4/\alpha}\theta-\cos^{8/\alpha}\theta-\sin^{8/\alpha}\theta-1}{4\cos^{4/\alpha}\theta}\\
&=\frac{2(\cos^{4/\alpha}\theta+\sin^{4/\alpha}\theta)-(\cos^{4/\alpha}\theta-\sin^{4/\alpha}\theta)^2-1}{4\cos^{4/\alpha}\theta}
\end{align*}
Thus we get $(w_1,w_2,w_3)=(v_1,v_2,v_3)$. \\
\eqref{A1a} Using the above calculation one can see that $(v_1,v_2,v_3)$ is indeed in the set $\mathcal{R}_{d,\alpha}$. 
\end{proof}

\begin{proof}[{\bf Proof of Proposition \ref{crs}}]

For $d=1$, the result follows from Lemma A
whereas for $d\geq2$ the result follows from Lemma B. 
\end{proof}

\subsection*{2. Comments on  geometric  optics  and  Fourier analytic approaches} \label{cop}
In this subsection, we compare geometric approach and Fourier analytic approach.  R\'emi Carles and his collaborators  in  \cite{ carles2012geometric, bhimani2020norm}  have established norm inflation with finite/infinite loss of  regularity for 
\begin{equation}\tag{A1}\label{dfe1}
i \partial_t \psi + \frac{1}{2} \Delta \psi = \pm |\psi|^{2\sigma}\psi, \quad \psi (0, x)= \psi_0
\end{equation}
  via geometric optics  approach in  $H^s$ for some $s<0.$   In this approach, the idea is to transfer the problem to a WKB analysis. This analysis seeks an approximate solution  to
\begin{equation}\tag{A2}\label{dfe2}
i \epsilon \partial_tu^{\epsilon} + \frac{\epsilon^2}{2} \Delta u^{\epsilon} =\pm  \epsilon |u^{\epsilon}|^{2\sigma} u^{\epsilon}
\end{equation}  
  like initial data.  And where the data is  given by  a superposition of highly oscillatory plane waves, i.e,  $u^{\epsilon}(0,x)=\sum_{j\in J_0\subset \mathbb N} \alpha_j(x) e^{i\kappa_j\cdot x/ \epsilon}.$ Here, wave vector $\kappa_j \in \mathbb R^d$ and amplitudes $\alpha_j\in \mathcal{S}(\R^d)$. Specifically, we seek, $$u^{\epsilon}(t,x)\approx \sum_{j\in J \supset J_0 } a_j(t,x)e^{i \phi_j(t,x)/\epsilon}.$$
  In order to do so, one proceeds as follows.
  \begin{itemize}
  \item[(i)]  We seek suitable scaling to  use  \textit{weakly non-linear  geometric optics}.  Specifically, one may find links between \eqref{dfe1} and \eqref{dfe2} via following relation: let $u^{\epsilon}(t,x)= \epsilon^{\frac{2-K}{2\sigma}} \psi(t,x) \ (\epsilon\to 0, 1\leq K<2)$ and for $\psi$ solution to \eqref{dfe1}, $u^{\epsilon}$ solves \eqref{dfe2}.
  \item[(ii)] Approximate  solution $u^{\epsilon}_{app}$ (of the exact solution $u^{\epsilon}$ of \eqref{dfe2}) given by geometric optics in the sense that $\|u^{\epsilon}-u^{\epsilon}_{app}\|_{L^{\infty}([0, T], L^2)}\to 0$ as $\epsilon\to 0.$ 
  \item[(iii)] As a first step  we  need to determine the characteristic  phases $\phi_j(t,x)\in \R.$ For the plane-wave initial data we are led to  eikonal  equation ($\partial_t \phi_j+ \frac{1}{2}|\nabla \phi_j|^2=0, \quad \phi_j(0,x)= \kappa_j\cdot x$).  The solution of which is explicitly given by 
\begin{equation*}
\phi_j(t,x)= \kappa_j\cdot x- \frac{t}{2}|\kappa_j|^2.
\end{equation*} 
 In general $J_0\subsetneq J,$ i.e, we usually need to take into account more phases  than we are given initially.
  The set of resonance  leading to the phase $\phi_j$ leading phase $\phi_j$ is therefore given by 
  \begin{align*}
 { \mathcal R}_j &= \left\{ (\kappa_\ell)_{1\le \ell\le 2\sigma+1}:\
                   \sum_{\ell=1}^{2\sigma+1} (-1)^{\ell+1} \kappa_\ell=\kappa_j,\
                   \sum_{\ell=1}^{2\sigma+1} (-1)^{\ell+1} |\kappa_\ell|^2=|\kappa_j|^2\right\}. 
\end{align*}
Continuing  the formal WKB approach, the $\mathcal{O}(\epsilon^1)$ term yields,  a system of transport equations.  This enables to choose initial profiles so that $\partial_ta_0\neq 0$ though $a_0$ is zero at initial time. This modes becomes instantaneously non-trivial.
  \item[(iv)] Establish a link between the norms of $\psi$ and $u^{\epsilon}_{app}.$ Finally,  steps (i) to (iii) yield that  $u_0$ and $u-u_{app}$ remains small in $H^s$ for $s<0$ but  $u_{app}$ is unbounded in $H^{\sigma}-$norm for all $\sigma\in \mathbb R$. This essentially  leads to the NI with infinite loss of regularity.
  \end{itemize}
(For some  $s>0$, several results for  norm inflation for NLS   are established in  \cite{carles2007geometric, thomann2008instabilities, alazard2009loss} via supercritical WKB analysis.)

Next, we comment on the Fourier analytic approach of the present paper: We take the  initial data of the form \eqref{kid} which are inverse Fourier transform of  characteristic  functions, whereas in WKB analysis, the initial data are in Schwartz class. Since we take data on the Fourier side as the characteristic function, the analysis is relatively  simplified and  can be more general. For the general strategy  for this approach we refer to Section \ref{mr}, the comments at the beginning of Section \ref{Sni}, and Remark \ref{psilg}. (See also \cite[Section 2]{kishimoto2019remark}.)

Note that we have proved the NI with infinite loss of regularity in two cases (Theorems \ref{mt2}, \ref{mt3} \eqref{mt32}). The proof of Theorem \ref{mt2} does not use the existence of resonant vectors. For the proof Theorem \ref{mt3} \eqref{mt32}  we indeed use it:
We highlight that in order to prove Lemma \ref{d2'}, the existence of $\eta_1,\eta_2,\eta_3\in\rd\setminus\{0\}$ satisfying the resonant condition
\[
\eta_1+\eta_2=\eta_3,\quad |\eta_1|^\alpha+|\eta_2|^\alpha=|\eta_3|^\alpha
\] has been crucially used.
 This should   compare with above point (iii). 
These vectors are used to get large  output in low frequency ${\xi\in Q_1}$ and they appear as a consequence of the non-linearity and the Fourier multiplier operator (free Schr\"odinger operator $U(t)$ given by $\F U(t)f=e^{it|\xi|^\alpha}\F f$).
 The advantage of this approach is that  we could    treat (fractional)  Hartree equation and NLS with power type nonlinearity simultaneously. 
 \\\\
 
\noindent 
\textbf{Acknowledgement}: Both authors are  thankful to Prof. K. Sandeep   for the encouragement  during this project and his  thoughtful suggestions on the preliminary draft of this paper.   Both authors are grateful to  Prof. R\'emi Carles for various thought comments and pointing out the connection between geometric and Fourier analytic approach. 
D.G. B is  thankful to DST-INSPIRE (DST/INSPIRE/04/2016/001507) for the research grant.  Both authors are  thankful to TIFR CAM for the  excellent  research facilities. Both the authors are grateful to an anonymous referee for carefully proof reading \&  pointing out  a flaw  in the argument in Lemma \ref{d0} in  previous version.\\

\noindent
\textbf{Declaration}: The authors have no conflict of interest.

\bibliographystyle{siam}
\bibliography{inflation}

\end{document}